\documentclass[12pt, reqno]{amsart}
\usepackage{amsfonts,latexsym, enumerate}
\usepackage{amsmath}
\usepackage{amscd}
\usepackage{float,amsmath,amssymb,mathrsfs,bm,multirow,graphics}
\usepackage[dvips]{graphicx}
\usepackage[percent]{overpic}
\usepackage[pdftex]{color}
\usepackage[all]{xy}

\newcommand{\nequation}{\setcounter{equation}{0}}
\renewcommand{\theequation}{\mbox{\arabic{section}.\arabic{equation}}}
\newcommand{\R}{{\Bbb R}}

\newcommand{\C}{{\Bbb C}}
\newcommand{\Z}{{\Bbb Z}}

\newcommand{\proofbegin}{\noindent{\it Proof.\quad}}
\newcommand{\proofend}{\hfill$\Box$\bigskip}

\newcommand{\I}{\mathbb{I}}

\DeclareMathOperator{\Ind}{Ind}
\DeclareMathOperator{\adj}{adj}

\newtheorem{theorem}{Theorem}[section]
\newtheorem{proposition}[theorem]{Proposition}

\newtheorem{lemma}[theorem]{Lemma}

\newtheorem{remark}[theorem]{Remark}
\newtheorem{example}[theorem]{Example}
\newtheorem{figuretext}{Figure}

\input epsf
\date{\today}
\title[Matrix Riemann-Hilbert problems]
{Matrix Riemann-Hilbert problems with jumps across Carleson contours}

\author{Jonatan Lenells}
\address{Department of Mathematics and Center for Astrophysics, Space Physics \& Engineering Research, Baylor University, One Bear Place \#97328, Waco, TX 76798, USA.}
\email{Jonatan\_Lenells@baylor.edu}

\begin{document}

\begin{abstract} \noindent
We develop a theory of $n \times n$-matrix Riemann-Hilbert problems for a class of jump contours and jump matrices of low regularity. Our basic assumption is that the contour $\Gamma$ is a finite union of simple closed Carleson curves in the Riemann sphere. In particular, contours with cusps, corners, and nontransversal intersections are allowed. We introduce a notion of $L^p$-Riemann-Hilbert problem and establish basic uniqueness results and a vanishing lemma. We also investigate the implications of Fredholmness for the unique solvability and prove a theorem on contour deformation.
\end{abstract}

\maketitle




\section{Introduction}\nequation
A Riemann-Hilbert (RH) problem consists of finding a sectionally analytic function with prescribed jumps across some given contour in the complex plane. In its simplest formulation, the problem involves a smooth simple closed contour $\Gamma$ dividing the complex plane into an interior domain $D_+$ and an exterior domain $D_-$, as well as a smooth `jump matrix' $v(z)$ defined for $z \in \Gamma$. 
The problem consists of finding an $n\times n$-matrix-valued function $m(z)$ which is analytic in $D_+ \cup D_-$ and whose boundary values $m_+$ and $m_-$ from the left and right sides of $\Gamma$ exist, are continuous, and satisfy the jump condition $m_+ = m_- v$ on $\Gamma$. Uniqueness is ensured by requiring that $m$ approaches the identity matrix at infinity. 

The theory of scalar RH problems is well-developed in the classical set-up in the complex plane \cite{M1992} as well as for problems on Riemann surfaces \cite{R1988, Z1971}. Constructive existence and uniqueness results are available, at least within classes of H\"older continuous functions \cite{AF2003, M1992, Z1971}. We refer to the monograph \cite{KKP1998} for more recent developments and further references in the case of less regular solutions.

The theory of matrix RH problems is substantially more complicated than the scalar theory. Only very special classes of problems (such as problems with a rational jump matrix, see Chapter I of \cite{CG1981}) can be solved explicitly. Uniqueness can often be established by means of Liouville's theorem, but existence results are rare and usually rely on the presence of some special symmetry, see \cite{AF2003, D1999}.

Matrix RH problems are essential in the analysis of integrable systems, orthogonal polynomials, and random matrices. The RH approach is particularly powerful when it comes to determining asymptotics. Indeed, the asymptotic behavior of solutions of many RH problems can be efficiently determined by means of the nonlinear steepest descent method introduced by Deift and Zhou \cite{DZ1993}, building on earlier work of Its \cite{I1981} and Manakov \cite{M1974}. This method and generalizations thereof have been instrumental in several recent advances in random matrix theory and in the analysis of large-time asymptotics of solutions of integrable PDE's \cite{D1999, DVZ1997, FI1996, FIKN2006, K2008, KT2012}.

The classical formulation of a RH problem, which involves a piecewise smooth contour $\Gamma$ and a smooth (or at least H\"older continuous) jump matrix $v$, is sufficient for many applications. However, in order to obtain a more convenient setting for the application of functional analytic techniques, it is essential to extend the formulation of a RH problem to the $L^p$-setting. Deift and Zhou and others \cite{D1999, DZ2002a, DZ2002b, FIKN2006, Z1989} have extended the definition of a RH problem to the case where the jump matrix $v$ and its inverse $v^{-1}$ belong to appropriate Lebesgue spaces, and the contour $\Gamma$ is a finite union of closed simple smooth curves in the Riemann sphere with a finite number of transversal intersection points. In particular, the relationship between the unique solvability of a RH-problem and the Fredholmness of a certain associated singular operator was explained in \cite{Z1989}. 

Our goal in this paper is to develop a theory of matrix RH problems for a class of jump contours of very low regularity. Our basic assumption is that the contour $\Gamma$ is a finite union of closed Carleson curves in the Riemann sphere. In particular, contours with cusps, corners, and nontransversal intersections are allowed. We introduce a notion of $L^p$-Riemann-Hilbert problem for this class of contours and establish basic uniqueness results and a vanishing lemma. We also investigate the implications of Fredholmness for the unique solvability and prove a theorem on contour deformation.

\begin{figure}
\begin{center}
 \begin{overpic}[width=.45\textwidth]{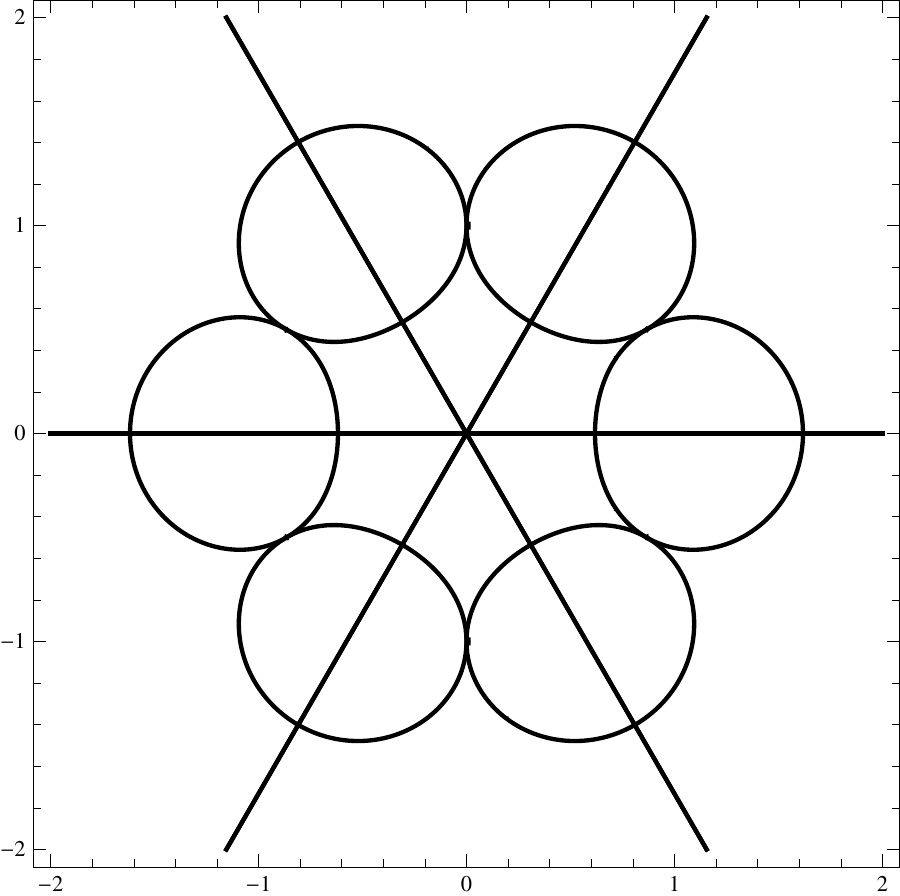}
    \end{overpic}
    \qquad \qquad
     \begin{figuretext}\label{Dns.pdf}
       A jump contour with nontransversal intersections that arises in the analysis of the Degasperis-Procesi equation on the half-line. 
         \end{figuretext}
     \end{center}
\end{figure}

We emphasize that RH problems with contours involving nontransversal intersections are important in applications to integrable evolution equations. For example, the analysis of the Degasperis-Procesi equation on the half-line naturally leads to a RH problem with the jump contour displayed in Figure \ref{Dns.pdf}, see \cite{L2013}. The results of this paper can be used to rigorously analyze the long-time asymptotics of the solutions of this equation via the nonlinear steepest descent method.

The formulation of a successful theory of RH problems is intricately linked to the boundedness of the Cauchy singular operator $\mathcal{S}_\Gamma$ defined in equation (\ref{cauchysingulardef}) below. Indeed, this operator is the key ingredient in the Sokhotski-Plemelj formulas for the boundary values of an analytic function. Since it has been proved in recent years that $\mathcal{S}_\Gamma$ is a bounded operator on $L^p(\Gamma)$ if and only if $\Gamma$ is Carleson (cf. \cite{BK1997}), it is natural to expect that the class of Carleson contours is the most general class of contours for which a clean RH theory exists. This is the reason for considering Carleson jump contours in this paper.

A second reason for writing this paper is to make accessible detailed and rigorous proofs of several basic results on matrix RH problem. Many of these results are well-known to the experts (at least if the contour $\Gamma$ is sufficiently well-behaved), but their proofs are scattered or absent in the literature. It turns out that the basic results can be proved in the more general setting of Carleson jump contours with little extra effort. 

In Section \ref{prelsec}, we summarize several properties of Smirnoff classes and Cauchy integrals over rectifiable Jordan curves.
In Section \ref{carlesonsec}, we introduce the notion of a Carleson jump contour as well as a number of function spaces which turn out to be convenient when dealing with contours passing through infinity. 
In Section \ref{cauchysec}, we establish several properties of Cauchy integrals over general Carleson jump contours. 
In Section \ref{rhsec}, we introduce a notion of $L^p$-Riemann-Hilbert problem for a general Carleson jump contour and develop the basics of a theory for these problems.

\section{Preliminaries}\label{prelsec}\nequation
A subset $\Gamma \subset \C$ is an {\it arc} if it is homeomorphic to a connected subset $I$ of the real line which contains at least two distinct points. If $\varphi:I \to \Gamma$ is a homeomorphism onto an arc and $(a,b) \subset I$ is the interior of $I$ with $a \in \R \cup \{-\infty\}$ and $b \in \R \cup \{\infty\}$, then $\lim_{t \to a^+} \varphi(t)$ and $\lim_{t \to b^-} \varphi(t)$ are referred to as {\it endpoints} of $\Gamma$ whenever the limits exist and are finite. An arc may have two, one, or no endpoints. An arc that does not contain its endpoints is an {\it open} arc. If $I = [a,b]$ is a closed interval, the {\it length} $|\Gamma|$ of $\Gamma$ is defined by
$$|\Gamma| = \sup \sum_{i=1}^n |\varphi(t_i) - \varphi(t_{i-1})|$$
where the supremum is over all partitions $a = t_0 < t_1 < \cdots < t_n = b$ of $[a,b]$. If $I$ is not closed, the length of $\Gamma$ is defined as the supremum of $|\varphi([a,b])|$ as $[a,b]$ ranges over all closed subintervals of $I$. The arc $\Gamma$ is {\it rectifiable} if its length is finite.
A subset $\Gamma \subset \C$ is a {\it composed curve} if it is connected and may be represented as the union of finitely many arcs each pair of which have at most endpoints in common. 
A composed curve is {\it oriented} if it can be represented as the union of finitely many oriented arcs each pair of which have at most endpoints in common. 
A subset $\Gamma \subset \C$ is a {\it Jordan curve} if it is homeomorphic to the unit circle $S^1$.

Let $\Gamma \subset \C$ be a composed curve. If $z \in \C$, $r \in (0, \infty)$, and $D(z,r)$ denotes the open disk of radius $r$ centered at $z$, then $\Gamma \cap D(z, r)$ is an at most countable union of arcs. If all of these arcs are rectifiable and the sum of their lengths is finite, we say that $\Gamma \cap D(z, r)$ is rectifiable. $\Gamma$ is {\it locally rectifiable} if $\Gamma \cap D(z, r)$ is rectifiable for every $z \in \Gamma$ and every $r \in (0,\infty)$.
A composed curve $\Gamma$ is locally rectifiable if and only if $\Gamma \cap D(0, r)$ is rectifiable for every $r \in (0,\infty)$.

\subsection{Carleson curves}
Let $\Gamma \subset \C$ be a locally rectifiable composed curve. We equip $\Gamma$ with Lebesgue length measure and denote the measure of a measurable subset $\gamma \subset \Gamma$ by $|\gamma|$; see e.g. Chapter 1 of \cite{P1956} for a detailed definition. We say that $\Gamma$ is {\it Carleson} (or a {\it Carleson curve}) if
\begin{align}\label{carlesondef}
 \sup_{z \in \Gamma} \sup_{r > 0} \frac{|\Gamma \cap D(z, r)|}{r} < \infty.
\end{align}
The condition (\ref{carlesondef}) is equivalent to the condition
\begin{align}\label{carlesondef2}
  \sup_{z \in \C} \sup_{r > 0} \frac{|\Gamma \cap D(z, r)|}{r} < \infty.
\end{align}
Moreover, $\Gamma$ is Carleson if and only if each of its finite number of arcs is Carleson.
We refer to Chapter 1 of \cite{BK1997} for more information on Carleson curves.

\subsection{Cauchy singular operator}
Let $\Gamma$ be a composed locally rectifiable curve. 
Let $C_0^\infty(\Gamma)$ denote the set of all restrictions of smooth functions $f:\R^2 \to \C$ of compact support to $\Gamma$. 
A measurable function $w:\Gamma \to [0,\infty]$ is a {\it weight} on $\Gamma$ if the preimage $w^{-1}(\{0,\infty\})$ has measure zero. The weighted Lebesgue space $L^p(\Gamma, w)$, $p \in [1, \infty)$, is defined as the space of measurable functions such that
$$\|f\|_{L^p(\Gamma, w)} := \bigg(\int_\Gamma |f(z)|^p w(z)^p |dz|\bigg)^{1/p} < \infty.$$
Equipped with the norm $\|\cdot \|_{L^p(\Gamma, w)}$, $L^p(\Gamma, w)$ is a complete Banach space. 

Given a function $h$ defined on $\Gamma$, we define $(\mathcal{S}_\Gamma h)(z)$ for $z \in \Gamma$ by
\begin{align}\label{cauchysingulardef}
(\mathcal{S}_\Gamma h)(z) = \lim_{\epsilon \to 0} \frac{1}{\pi i} \int_{\Gamma \setminus D(z, \epsilon)} \frac{h(z')}{z' - z} dz',
\end{align}
whenever the limit exists. 
If $h \in C_0^\infty(\Gamma)$, then $(\mathcal{S}_\Gamma h)(z)$ exists for almost all $z \in \Gamma$ (see Theorem 4.14 of \cite{BK1997}). If $w$ is a weight on $\Gamma$, we say that the Cauchy singular operator $\mathcal{S}_\Gamma$ {\it generates a bounded operator on $L^p(\Gamma, w)$} if $C_0^\infty(\Gamma)$ is dense in $L^p(\Gamma, w)$ and 
$$\|\mathcal{S}_\Gamma f\|_{L^p(\Gamma, w)} < M \|f\|_{L^p(\Gamma,w)} \quad \text{for all} \quad f \in C_0^\infty(\Gamma)$$
with some constant $M > 0$ independent of $f$; in that case there exists a unique bounded operator $\tilde{\mathcal{S}}_\Gamma$ on $L^p(\Gamma, w)$ such that $\tilde{\mathcal{S}}_\Gamma f = \mathcal{S}_\Gamma f$ for all $f \in C_0^\infty(\Gamma)$.

It was realized in the early eighties that the Carleson condition is the essential condition for ascertaining boundedness of $\mathcal{S}_\Gamma$ in $L^p$-spaces \cite{D1984}. More precisely, if $1 < p < \infty$ and $\Gamma$ is a composed locally rectifiable curve, then $\mathcal{S}_\Gamma$ generates a bounded operator on $L^p(\Gamma)$ if and only if $\Gamma$ is Carleson. In dealing with unbounded contours, we will need a more general version of this result valid for weighted $L^p$-spaces.

Let $p \in (1, \infty)$. We define $A_p(\Gamma)$ as the set of weights $w \in L_{\text{loc}}^p(\Gamma)$ such that $1/w \in L_{\text{loc}}^q(\Gamma)$ and
\begin{align}\label{muckenhouptcondition}
\sup_{z \in \Gamma} \sup_{r >0} \bigg(\frac{1}{r} \int_{\Gamma \cap D(z, r)} w(z')^p |dz'|\bigg)^{1/p}\bigg(\frac{1}{r} \int_{\Gamma \cap D(z, r)} w(z')^{-q} |dz'|\bigg)^{1/q} < \infty,
\end{align}
where $q \in (1,\infty)$ is defined by $1/p + 1/q = 1$.
Elements in $\cup_{1<p<\infty} A_p(\Gamma)$ are referred to as {\it Muckenhoupt weights} on $\Gamma$. If $\Gamma$ is Carleson, then constant weights belong to $A_p(\Gamma)$. If $A_p(\Gamma)$ is nonempty, then $\Gamma$ is Carleson.

\begin{theorem}\label{Sth}
Let $1 < p < \infty$ and let $\Gamma$ be a composed locally rectifiable curve. 
Let $w$ be a weight on $\Gamma$. Then $\mathcal{S}_\Gamma$ generates a bounded operator $\tilde{\mathcal{S}}_\Gamma$ on $L^p(\Gamma, w)$ if and only if $\Gamma$ is Carleson and $w \in A_p(\Gamma)$. Moreover, if $f \in L^p(\Gamma, w)$ and  $\mathcal{S}_\Gamma$ generates a bounded operator on $L^p(\Gamma, w)$, then the limit in (\ref{cauchysingulardef}) exists and $(\mathcal{S}_\Gamma f)(z) = (\tilde{\mathcal{S}}_\Gamma f)(z)$ for a.e. $z \in\Gamma$.
\end{theorem}
\proofbegin
See Theorem 4.15 and Remark 5.23 of \cite{BK1997}.
\proofend

\subsection{Smirnoff classes}\label{smirnoffsubsec}
Let $\Gamma \subset \C$ be a rectifiable Jordan curve oriented counterclockwise. Let $\hat{\C} = \C\cup \infty$ denote the Riemann sphere and let $D_+$ and $D_-$ be the two components of $\hat{\C} \setminus \Gamma$. Assuming that $\infty \in D_-$, we refer to $D_+$ and $D_-$ as the interior and exterior components respectively.
Let $1 \leq p < \infty$. A function $f$ analytic in $D_+$ belongs to the {\it Smirnoff class $E^p(D_+)$} if there exists a sequence of rectifiable Jordan curves $\{C_n\}_1^\infty$ in $D_+$, tending to the boundary in the sense that $C_n$ eventually surrounds each compact subdomain of $D_+$, such that
\begin{align}\label{Epsup}
\sup_{n \geq 1} \int_{C_n} |f(z)|^p |dz| < \infty.
\end{align}
A function $f$ analytic in $D_-$ is said to be of class $E^p(D_-)$ if there exists a sequence of rectifiable Jordan curves $\{C_n\}_1^\infty$ in $D_-$, tending to the boundary $\Gamma$ in the sense that every compact subset of $D_-$ eventually lies outside $\Gamma_n$, such that (\ref{Epsup}) holds.
We let $\dot{E}^p(D_-)$ denote the subspace of $E^p(D_-)$ consisting of all functions $f \in E^p(D_-)$ that vanish at infinity.

\subsection{Basic results on Cauchy integrals}
Given a locally rectifiable composed contour $\Gamma \subset \C$ and a measurable functionÊ $h$ defined on $\Gamma$, we define the Cauchy integral $(\mathcal{C}h)(z)$ for $z \in \C \setminus \Gamma$ by
\begin{align}\label{cauchyintegraldef}
(\mathcal{C}h)(z) = \frac{1}{2\pi i} \int_\Gamma \frac{h(z')dz'}{z' - z},
\end{align}
whenever the integral converges. To avoid confusion, we will sometimes indicate the dependence of $\mathcal{C}$ on $\Gamma$ by writing $\mathcal{C}_\Gamma$Ê for $\mathcal{C}$.

In the next two propositions, we collect a number of properties of the Cauchy integral and its relation to the Smirnoff classes; we refer to Chapter 10 of \cite{D1970} and Chapter 6 of \cite{BK1997} for proofs. Given a Jordan curve $\Gamma \subset \C$, we let $D_+$ and $D_-$ denote the interior and exterior components of $\hat{\C} \setminus \Gamma$. 

\begin{theorem}\label{Jordanth1}
Let $\mathcal{C}$ denote the Cauchy integral operator defined in (\ref{cauchyintegraldef}).
\begin{enumerate}[$(a)$]
\item Let $1 \leq p < \infty$. Suppose $\Gamma \subset \C$ is a rectifiable Jordan curve.
If $f \in E^p(D_+)$, then the nontangential limits of $f(z)$ as $z$ approaches the boundary exist a.e. on $\Gamma$; if $f_+(z)$ denotes the boundary function, then $f_+ \in L^p(\Gamma)$ and
\begin{align}\label{CGammaf}
(\mathcal{C} f_+)(z) = \begin{cases} f(z), & z \in D_+, \\
0, & z \in D_-.
\end{cases}
\end{align}
If $f \in E^p(D_-)$, then the nontangential limits of $f(z)$ as $z$ approaches the boundary exist a.e. on $\Gamma$. If $f_-(z)$ denotes the boundary function, then $f_- \in L^p(\Gamma)$ and
\begin{align}\label{CGammaminus}
(\mathcal{C} f_-)(z) = \begin{cases}  f(\infty), & z \in D_+, \\
f(\infty) - f(z), & z \in D_-.
\end{cases}
\end{align}

\item Let $1 < p < \infty$. Suppose $\Gamma$ is a Carleson Jordan curve. Then the Cauchy singular operator $\mathcal{S}_\Gamma:L^p(\Gamma) \to L^p(\Gamma)$ defined in (\ref{cauchysingulardef}) satisfies $\mathcal{S}_\Gamma^2 = I$. Moreover,
if $h \in L^p(\Gamma)$, then 
$$(\mathcal{C}h)|_{D_+} \in E^p(D_+), \qquad (\mathcal{C}h)|_{D_-} \in \dot{E}^p(D_-).$$
\end{enumerate}
\end{theorem}

Theorem \ref{Jordanth1} implies that if $\Gamma$ is a Carleson Jordan curve and $h \in L^p(\Gamma)$ for some $1 < p < \infty$, then the left and right nontangential boundary values of $\mathcal{C}h$, which we denote by $\mathcal{C}_+ h$ and $\mathcal{C}_- h$, lie in $L^p(\Gamma)$. This allows us to view $\mathcal{C}_\pm$ as linear operators $\mathcal{C}_\pm:h \mapsto \mathcal{C}_\pm h$ on $L^p(\Gamma)$.

\begin{theorem}\label{Jordanth2}
Let $1 < p < \infty$ and let $\Gamma \subset \C$ be a Carleson Jordan curve.
Then $\mathcal{C}_\pm$ are bounded operators on $L^p(\Gamma)$ with the following properties:

\begin{itemize}
\item The Sokhotski-Plemelj formulas 
$$\mathcal{C}_+ = \frac{1}{2}(I + \mathcal{S}_\Gamma), \qquad \mathcal{C}_- = \frac{1}{2}(-I + \mathcal{S}_\Gamma),$$
are valid. 

\item $\mathcal{C}_\pm$ are complementary projections on $L^p(\Gamma)$ in the sense that
$$L^p(\Gamma) = \mathcal{C}_+L^p(\Gamma) \oplus \mathcal{C}_-L^p(\Gamma)$$ 
and
$$\mathcal{C}_+ - \mathcal{C}_- = \I, \qquad \mathcal{C}_+^2 = \mathcal{C}_+, \qquad \mathcal{C}_-^2 = -\mathcal{C}_-, \qquad \mathcal{C}_+\mathcal{C}_- = \mathcal{C}_-\mathcal{C}_+ = 0.$$

\item If $h = \mathcal{C}_+h - \mathcal{C}_-h \in L^p(\Gamma)$, then
$$(\mathcal{C}h)|_{D_+} = (\mathcal{C}\mathcal{C}_+h)|_{D_+} \in E^p(D_+), \qquad (\mathcal{C}h)|_{D_-} = -(\mathcal{C}\mathcal{C}_-h)|_{D_-} \in \dot{E}^p(D_-).$$

\item The map $h \mapsto (\mathcal{C}h)|_{D_+}$ is a bijection  $\mathcal{C}_+L^p(\Gamma) \to E^p(D_+)$ with inverse $f \mapsto f_+$. 

\item The map $h \mapsto (\mathcal{C}h)|_{D_-}$ is a bijection  $\mathcal{C}_-L^p(\Gamma) \to \dot{E}^p(D_-)$ with inverse $f \mapsto -f_-$.
\end{itemize}
\end{theorem}

\section{Carleson jump contours}\label{carlesonsec}\nequation
RH problems are conveniently formulated on the Riemann sphere $\hat{\C} = \C\cup \infty$.
In order to allow for jump contours passing through infinity, we introduce a class of curves $\mathcal{J}$, which in addition to the rectifiable Jordan curves considered in the previous section also includes unbounded contours. Recall that $\Gamma \subset \C$ is referred to as a Carleson curve if $\Gamma$ is a locally rectifiable composed curve satisfying (\ref{carlesondef}). We extend this notion to the Riemann sphere by calling a subset $\Gamma \subset \hat{\C}$ a Carleson curve if and only if $\Gamma$ is connected and $\Gamma \cap \C$ is a Carleson curve.

\subsection{The class $\mathcal{J}$}
Let $\mathcal{J}$ denote the collection of all subsets $\Gamma$ of the Riemann sphere $\hat{\C}$ such that $\Gamma$ is homeomorphic to the unit circle and $\Gamma$ is a Carleson curve. 
If $\infty \notin \Gamma$, then $\Gamma \in \mathcal{J}$ if and only if $\Gamma \subset \C$ is a Carleson Jordan curve. 
However, $\mathcal{J}$ also includes curves passing through infinity. 
In fact, the next proposition shows that $\mathcal{J}$ is invariant under the action of the group of linear fractional transformations. This shows that $\mathcal{J}$ is a natural extension of the family of Carleson Jordan curves in that it puts $\infty$ on an equal footing with the other points in the Riemann sphere.

\begin{proposition}\label{Jprop}
The family of all Carleson curves in $\hat{\C}$ is invariant under the action of the group of linear fractional transformations. In other words, if $\psi: \hat{\C} \to \hat{\C}$ is given by
\begin{align}\label{linearfracdef}
\psi(z) = \frac{az + b}{cz + d},
\end{align}
for some constants $a,b,c,d \in \C$ with $ad - bc \neq 0$, then $\Gamma \subset \hat{\C}$ is a Carleson curve if and only if $\psi(\Gamma) \subset \hat{\C}$ is a Carleson curve.
\end{proposition}
\proofbegin
See Appendix \ref{Japp}.
\proofend

\begin{remark}\upshape
If $\gamma:S^1 \to \hat{\C}$ is an injective $C^1$ map such that $|\gamma'(s)| \neq 0$ for all $s \in S^1$, then $\Gamma := \gamma(S^1)$ belongs to $\mathcal{J}$. Indeed, being a continuous bijection from a compact space onto a Hausdorff space, $\gamma$ is a homeomorphism $S^1 \to \Gamma$.
In view of Proposition \ref{Jprop}, we may assume that $\infty \notin \Gamma$. Then, since $S^1$ is compact, we may cover $S^1$ with a finite number of open sets $\{U_j\}_1^n$ such that the restriction of $\gamma$ to each $U_j$ is a $C^1$ graph; Proposition 1.1 of \cite{BK1997} now implies that $\Gamma$ is Carleson.
\end{remark}

\begin{remark}\upshape
The Carleson condition is essential in Proposition \ref{Jprop}. In fact, the family of composed locally rectifiable (but not necessarily Carleson) curves is {\it not} invariant under the action of the group of linear fractional transformations. Indeed, let $\Gamma = \{te^{-it^2} | 1 < t < \infty\}$ and $\psi(z) = z^{-1}$. Then $\Gamma$ is locally rectifiable, but $\psi(\Gamma)$ is not locally rectifiable because $\psi(\Gamma) = \{t^{-1}e^{it^2} | 1 < t < \infty\} \subset D(0,1)$ has infinite length:
$$\int_1^\infty \bigg|\frac{d}{dt} t^{-1}e^{it^2} \bigg| dt
= \int_1^\infty \sqrt{4 + \frac{1}{t^4}} dt = \infty.$$
This example does not contradict Proposition \ref{Jprop}. Indeed, the estimate 
$$\frac{|\Gamma \cap D(0, r)|}{r}
= \frac{1}{r}\int_1^r \bigg|\frac{d}{dt} te^{-it^2}\bigg| dt
= \frac{1}{r}\int_1^r \sqrt{1 + 4t^4} dt
> \frac{1}{r} \int_1^r 2t^2 dt = \frac{2(r^3 - 1)}{3r}$$
implies that $\frac{|\Gamma \cap D(0, r)|}{r}$ is unbounded as $r \to \infty$; hence $\Gamma$ is not Carleson.
\end{remark}

\subsection{Carleson jump contours}

We call a connected subset $\Gamma \subset \hat{\C}$ of the Riemann sphere a {\it Carleson jump contour} if it has the following properties:
\begin{enumerate}[$(a)$]
\item $\Gamma \cap \C$ is an oriented composed curve.

\item $\hat{\C} \setminus \Gamma$ is the union of two disjoint open sets $D_+$ and $D_-$ each of which has a finite number of simply connected components in $\hat{\C}$.

\item $\Gamma$ is the positively oriented boundary of $D_+$ and the negatively oriented boundary of $D_-$, i.e. $\Gamma = \partial D_+ = -\partial D_-$.

\item If $\{D_j^+\}_1^n$ and $\{D_j^-\}_1^m$ are the components of $D_+$ and $D_-$, then $\partial D_j^+ \in \mathcal{J}$ for $j = 1, \dots, n$, and $\partial D_j^- \in \mathcal{J}$ for $j = 1, \dots, m$.
\end{enumerate}

\begin{figure}
\bigskip\bigskip
\begin{center}
 \begin{overpic}[width=.6\textwidth]{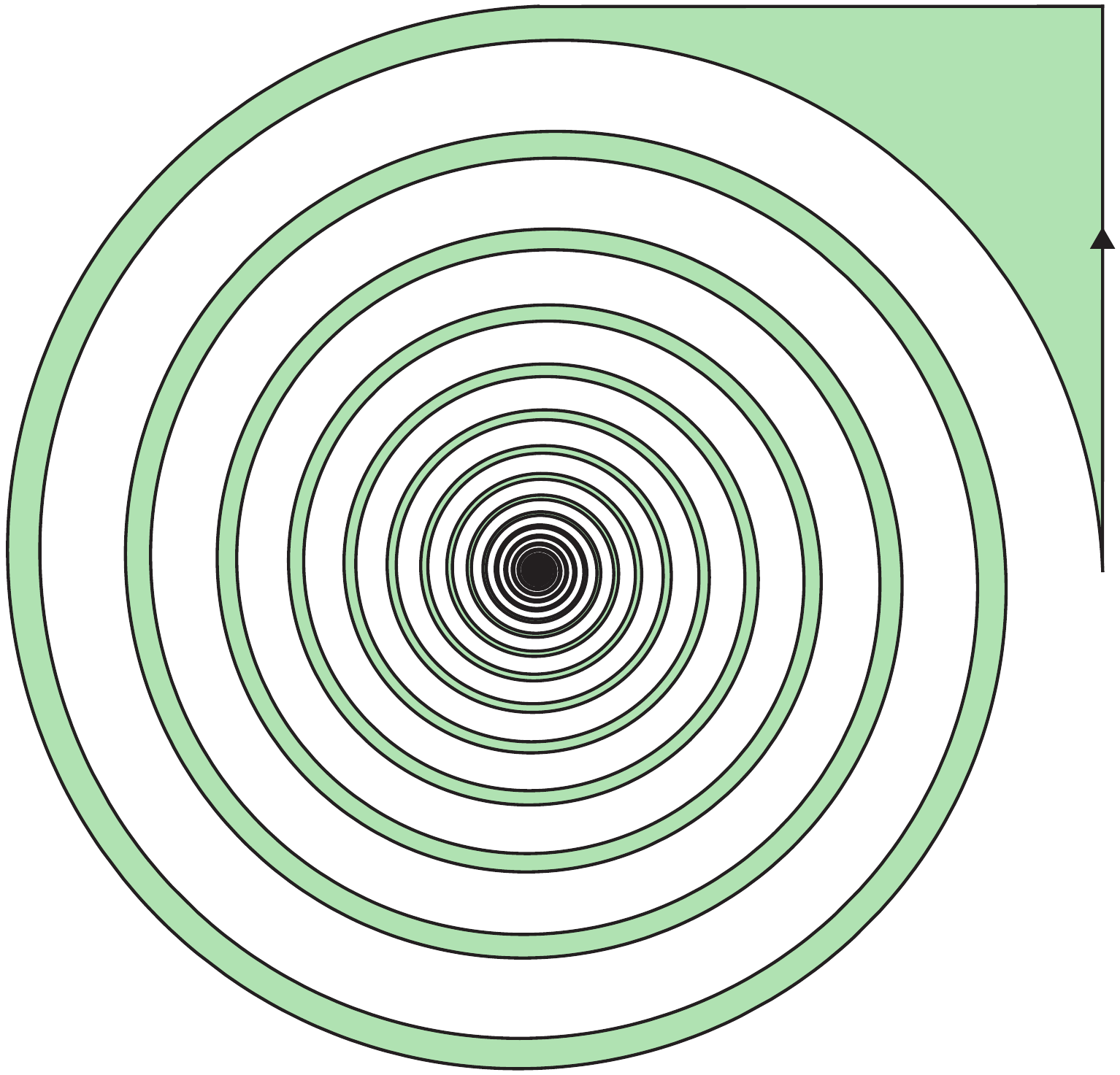}
 \put(87,85){$D_+$}
 \put(107,85){$D_-$}
 \put(101,73){$\Gamma$}
 \put(101,44){$1$}
 \put(100,97){$1+i$}
 \put(47,98){$i$}
 \end{overpic}
   \bigskip
   \begin{figuretext}\label{spiral.pdf}
    The Carleson jump contour defined in (\ref{spiraldef}). 
    \end{figuretext}
   \end{center}
\end{figure}

\begin{figure}
\begin{center}
 \begin{overpic}[width=.5\textwidth]{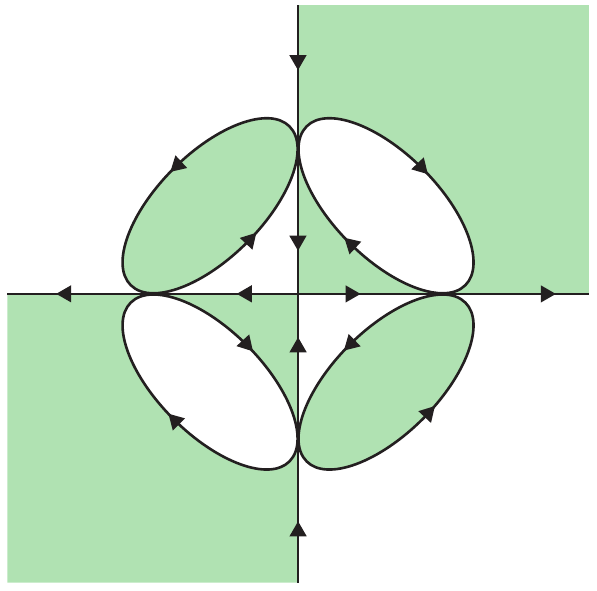}
 \put(80,80){$D_+$}
 \put(80,16){$D_-$}
 \put(102,48){$\Gamma$}
 \end{overpic}
   \bigskip
   \begin{figuretext}\label{jump1.pdf}
   An example of a Carleson jump contour.
      \end{figuretext}
   \end{center}
\end{figure}

\begin{figure}
\begin{center}
 \begin{overpic}[width=.65\textwidth]{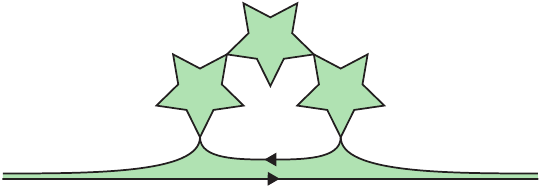}
 \put(47,27){$D_+$}
 \put(80,16){$D_-$}
 \put(49,-4){$\Gamma$}
 \end{overpic}
   \bigskip
   \begin{figuretext}\label{jump2.pdf}
   An example of a Carleson jump contour.
      \end{figuretext}
   \end{center}
\end{figure}

\begin{example}\upshape
Define the curve $\Gamma \subset \C$ by (see Figure \ref{spiral.pdf})
\begin{align}\label{spiraldef}
\Gamma = [1, 1+i] \cup [1+i, i] \cup \{ire^{-25 i \ln r} \, | \, 0 < r < 1\} \cup \{0\} \cup \{re^{-25 i \ln r} \, | \, 0 < r < 1\},
\end{align}
where $[a,b]$ denotes the straight line segment from $a$ to $b$. Then $\Gamma$ is a Carleson jump contour. Indeed, using that any logarithmic spiral $\{re^{-i\delta \ln r} \, | \, 0 < r < 1\}$ for $\delta \in \R$ is a Carleson arc (see Example 1.6 of \cite{BK1997}), it is straightforward to show that $\Gamma$ is a Carleson Jordan curve. 
Other examples of Carleson jump contours are displayed in Figures \ref{jump1.pdf} and \ref{jump2.pdf}.
\end{example}

Proposition \ref{Jprop} implies the following result.

\begin{proposition}\label{jumpprop}
  The family of Carleson jump contours is invariant under the action of the group of linear fractional transformations. In other words, if $\psi: \hat{\C} \to \hat{\C}$ is given by (\ref{linearfracdef})
for some constants $a,b,c,d \in \C$ with $ad - bc \neq 0$, then $\Gamma$ is a Carleson jump contour if and only if $\psi(\Gamma)$ is a Carleson jump contour.
\end{proposition}

Our goal is to establish generalizations of Theorems \ref{Jordanth1} and \ref{Jordanth2} which are valid in the case of a general Carleson jump contour $\Gamma$. These generalizations will be stated and proved in Section Ê\ref{cauchysec}; in the remainder of this section, we introduce a number of function spaces which will be needed for the formulation of these theorems. 

\subsection{Generalized Smirnoff classes}
In Subsection \ref{smirnoffsubsec}, we defined the Smirnoff class $E^p(D)$ for $D = D_+$ and $D = D_-$, where $D_+$ and $D_-$ are the domains interior and exterior to a rectifiable Jordan curve, respectively. We now extend this definition to allow for situations where $D$ is an arbitrary finite disjoint union of domains bounded by curves in $\mathcal{J}$.

Let $1 \leq p < \infty$. If $D$ is a subset of $\hat{\C}$ bounded by a curve in $\mathcal{J}$ which passes through infinity, we define $E^p(D)$ as the set of functions $f$ analytic in $D$ for which $f \circ \varphi^{-1} \in E^p(\varphi(D))$, where 
\begin{align}\label{varphidef}
  \varphi(z) = \frac{1}{z - z_0}
\end{align}  
and $z_0$ is any point in $\C \setminus \Gamma$. It is easy to see that $E^p(D)$ does not depend on the choice of $z_0$.
If $D$ is a subset of $\hat{\C}$ bounded by a curve in $\mathcal{J}$, we define $\dot{E}^p(D)$ as the subspace of $E^p(D)$ consisting of all functions $f \in E^p(D)$ such that $z f(z) \in E^p(D)$. If $D$ is bounded, then $\dot{E}^p(D) = E^p(D)$. If $\infty \in D$, then $\dot{E}^p(D)$ consists of the functions in $E^p(D)$ that vanish at infinity, so that the present definition of $\dot{E}^p(D)$ is consistent with the definition given in Subsection \ref{smirnoffsubsec}.  

If $D = D_1 \cup \cdots \cup D_n$ is the union of a finite number of disjoint subsets of $\hat{\C}$ each of which is bounded by a curve in $\mathcal{J}$, we define $E^p(D)$ and $\dot{E}^p(D)$ as the set of functions $f$ analytic in $D$ such that $f|_{D_j} \in E^p(D_j)$ and $f|_{D_j} \in \dot{E}^p(D_j)$ for each $j$, respectively.

\subsection{Properties of $E^p(D)$ and $\dot{E}^p(D)$}
Our definitions of the generalized Smirnoff classes $E^p(D)$ and $\dot{E}^p(D)$ are designed in such a way that these classes possess convenient transformation properties under the action of the group of linear fractional transformations.

\begin{proposition}\label{EpDinvariantprop}
Let $1 \leq p < \infty$. Let $\Gamma$ be a Carleson jump contour and let $D$ be the union of any number of components of $\hat{\C} \setminus \Gamma$. Let $\psi(z)$ be an arbitrary linear fractional transformation of the form (\ref{linearfracdef}) with $ad - bc \neq 0$. 
\begin{enumerate}[$(a)$]
\item $f \in E^p(D)$ if and only if $f \circ \psi^{-1} \in E^p(\psi(D))$. 

\item $f \in \dot{E}^p(D)$ if and only if $\Psi f \in \dot{E}^p(\psi(D))$ where $(\Psi f)(w) = (cw - a)^{-1} f(\psi^{-1}(w))$.
\end{enumerate}
\end{proposition}
\proofbegin
Without loss of generality, we may assume that $D$ is one of the components ofÊ $\hat{\C} \setminus \Gamma$ where $\Gamma \in \mathcal{J}$.

$(a)$
We will prove that $f \in E^p(D)$ if and only if $f \circ \psi^{-1} \in E^p(\psi(D))$ whenever $\psi(D)$ is bounded and either $(i)$ $\infty \in D$, $(ii)$ $D$ is bounded, or $(iii)$ $\infty \in \Gamma$.
Since the linear fractional transformations form a group, this will prove $(a)$.
 
{\it Case $(i)$.} Suppose $\psi(D)$ is bounded, $\infty \in D$, and $f \in E^p(D)$.
By the definition of $E^p(D)$ in Subsection \ref{smirnoffsubsec}, there exists a sequence of rectifiable Jordan curves $\{C_n\}_1^\infty$ in $D$ tending to the boundary $\Gamma$ such that (\ref{Epsup}) holds. It follows that $\{\psi(C_n)\}_1^\infty$ is a sequence of rectifiable Jordan curves in $\psi(D)$ tending to $\psi(\Gamma)$ and
\begin{align}\label{subsup}
\sup_{n \geq 1} \int_{\psi(C_n)} |f(\psi^{-1}(w))|^p |dw| 
= \sup_{n \geq 1} \int_{C_n} |f(z)|^p |\psi'(z)| |dz|.
\end{align}
If $c = 0$, then $\psi'(z)$ is a finite constant. If $c \neq 0$, then our assumption that $\psi(D)$ is bounded implies that the point $z = -d/c$ does not belong to $\bar{D}$; hence the function
$$\psi'(z) = \frac{ad - bc}{(cz + d)^2}$$
is bounded on $D$. It follows that the right-hand side of (\ref{subsup}) is finite. Thus $f \circ \psi^{-1} \in E^p(\psi(D))$. Conversely, if $\psi(D)$ is bounded, $\infty \in D$, and $f \circ \psi^{-1} \in E^p(\psi(D))$, then a similar argument shows that $f \in E^p(D)$.

{\it Case $(ii)$.}  Suppose both $\psi(D)$ and $D$ are bounded. An argument similar to that used in Case $(ii)$ shows that $f \in E^p(D)$ if and only if $f \circ \psi^{-1} \in E^p(\psi(D))$.

{\it Case $(iii)$.}  Suppose $\psi(D)$ is bounded and $\infty \in \Gamma$. 
Let $z_0 \in \C \setminus (D \cup \Gamma)$. By the definition of $E^p(D)$, $f \in E^p(D)$ if and only if $f \circ \varphi^{-1} \in E^p(\varphi(D))$, where $\varphi(z) = 1/(z - z_0)$. But $\psi \circ \varphi^{-1}$ is a linear fractional transformation mapping the bounded domain $\varphi(D)$ onto the bounded domainÊ $\psi(D)$; hence Case $(ii)$ implies that $f \circ \varphi^{-1} \in E^p(\varphi(D))$ if and only if
$f \circ \varphi^{-1} \circ (\psi \circ \varphi^{-1})^{-1} = f \circ \psi^{-1}$ belongs to $E^p(\psi(D))$.
This completes the proof of $(a)$.
 
$(b)$ This part is a consequence of $(a)$ and the definitions. Indeed, suppose first that $c \neq 0$. By definition of $\dot{E}^p(D)$, $f \in \dot{E}^p(D)$ if and only if $f(z), zf(z) \in E^p(D)$. Since $E^p(D)$ is a linear space, this is the case if and only if $f(z), (cz + d)f(z) \in E^p(D)$. Using $(a)$ and the fact that $(cz + d) = \frac{bc - ad}{cw - a}$ when $w = \psi(z)$, the latter condition is equivalent to $f(\psi^{-1}(w)), \frac{bc - ad}{cw - a}f(\psi^{-1}(w)) \in E^p(\psi(D))$. Using that $E^p(D)$ is a linear space again, this holds if and only if $(cw - a)^{-1} f(\psi^{-1}(w)) \in \dot{E}^p(\psi(D))$. The proof when $c = 0$ is similar. This proves $(b)$.
\proofend

\begin{lemma}\label{EpCnlemma}
Let $1 \leq p < \infty$. Let $D$ be a subset of $\hat{\C}$ bounded by a curve $\Gamma \in \mathcal{J}$. Let $z_0 \in \C \setminus \bar{D}$ and let $f:D \to \C$ be an analytic function.
\begin{enumerate}[$(a)$]
\item If there exist curves $\{C_n\}_1^\infty \subset \mathcal{J}$ in $D$, tending to $\Gamma$ in the sense that $C_n$ eventually surrounds each compact subset of $D \subset \hat{\C}$, such that
\begin{align}\label{Epsupz0}
\sup_{n \geq 1} \int_{C_n} |z - z_0|^{-2} |f(z)|^p |dz| < \infty,
\end{align}
then $f \in E^p(D)$.

\item If there exist curves $\{C_n\}_1^\infty \subset \mathcal{J}$ in $D$, tending to $\Gamma$ in the sense that $C_n$ eventually surrounds each compact subset of $D \subset \hat{\C}$, such that
\begin{align}\label{Epdotsupz0}
\sup_{n \geq 1} \int_{C_n} |z - z_0|^{p-2} |f(z)|^p |dz| < \infty,
\end{align}
then $f \in \dot{E}^p(D)$.
\end{enumerate}
\end{lemma}
\proofbegin
$(a)$ Let $\varphi(z) = \frac{1}{z-z_0}$. If $\{C_n\}_1^\infty \subset \mathcal{J}$ are curves in $D$ tending to $\Gamma$ such that (\ref{Epsupz0}) holds, then $\{\varphi(C_n)\}_1^\infty$ are curves in $\varphi(D)$ tending to $\varphi(\Gamma)$ and
$$\sup_{n \geq 1} \int_{\varphi(C_n)} |f(\varphi^{-1}(w))|^p |dw|
= \sup_{n \geq 1} \int_{C_n} |z - z_0|^{-2} |f(z)|^p |dz| < \infty.$$
Since $\varphi(C_n) \in \mathcal{J}$ by Proposition \ref{Jprop}, each $\varphi(C_n)$ is a rectifiable Jordan curve. Hence $f \circ \varphi^{-1} \in E^p(\varphi(D))$, that is, $f \in E^p(D)$.

$(b)$ If $\{C_n\}_1^\infty \subset \mathcal{J}$ are curves in $D$ tending to $\Gamma$ such that (\ref{Epdotsupz0}) holds, then $\{\varphi(C_n)\}_1^\infty$ are rectifiable Jordan curves in $\varphi(D)$ tending to $\varphi(\Gamma)$ and
$$\sup_{n \geq 1} \int_{\varphi(C_n)} |w^{-1} f(\varphi^{-1}(w))|^p |dw|
= \sup_{n \geq 1} \int_{C_n} |z - z_0|^{p-2} |f(z)|^p |dz| < \infty.$$
Hence $w^{-1} f(\varphi^{-1}(w)) \in E^p(\varphi(D))$; thus $f \in \dot{E}^p(D)$ by Proposition \ref{EpDinvariantprop}.
\proofend

\begin{lemma}\label{Epdotlemma}
Let $D$ be a subset of $\hat{\C}$ bounded by a curve in $\mathcal{J}$.
\begin{enumerate}[$(a)$]
\item $\dot{E}^p(D) \subset \dot{E}^r(D)$ whenever $1 \leq r \leq p < \infty$.
\item Suppose $p,q,r \in [1, \infty)$ satisfy $1/p + 1/q = 1/r$. If $f \in \dot{E}^p(D)$ and $g \in \dot{E}^q(D)$, then the functions $zf(z)g(z)$ and $fg$ belong to $\dot{E}^r(D)$.
\end{enumerate}
\end{lemma}
\proofbegin
$(a)$ The result is immediate from the definitions if $\infty \notin \Gamma$. Thus suppose $\infty \in \Gamma$. Let $z_0 \in \C \setminus \Gamma$ and define $\varphi$ as in (\ref{varphidef}). If $f \in \dot{E}^p(D)$, then Proposition \ref{EpDinvariantprop} implies that $w^{-1}f(\varphi^{-1}(w)) \in \dot{E}^p(\varphi(D))$; since $\dot{E}^p(\varphi(D)) \subset \dot{E}^r(\varphi(D))$, another application of Proposition \ref{EpDinvariantprop} shows thatÊ $f \in \dot{E}^r(D)$. This proves $(a)$. 

$(b)$ Suppose $p,q,r \in [1, \infty)$ satisfy $1/p + 1/q = 1/r$. Let $f \in \dot{E}^p(D)$ and $g \in \dot{E}^q(D)$. 
We first suppose $D$ is bounded. Then there exist sequences of rectifiable Jordan curves $\{A_n\}_1^\infty$ and $\{B_n\}_1^\infty$ in $D$ tending to the boundary of $D$ such that
\begin{align*}
\sup_{n \geq 1} \|f\|_{L^p(A_n)} < \infty, \qquad
\sup_{n \geq 1} \|g\|_{L^q(B_n)} < \infty.
\end{align*} 
Without loss of generality, we may assume that $A_n = B_n = C_n$ for all $n \geq 1$ where $\{C_n\}_1^\infty$ are level curves of an arbitrary conformal map of the unit disk onto $D$ (see Theorem 10.1 in \cite{D1970}). 
Then, by H\"older's inequality,
$$\sup_{n \geq 1} \|fg\|_{L^r(C_n)} \leq \sup_{n \geq 1} \|f\|_{L^p(C_n)} \|g\|_{L^q(C_n)} < \infty.$$
Hence $fg \in E^r(D) = \dot{E}^r(D)$. This proves $(b)$ when $D$ is bounded.

If $D$ is unbounded, then pick $z_0 \in \C \setminus (D \cup \Gamma)$ and let $\varphi(z) = 1/(z - z_0)$. By Proposition \ref{EpDinvariantprop}, $w^{-1} f(\varphi^{-1}(w)) \in \dot{E}^p(\varphi(D))$ and $w^{-1} g(\varphi^{-1}(w)) \in \dot{E}^q(\varphi(D))$. Hence, by the preceding paragraph, 
$w^{-2}(fg)(\varphi^{-1}(w)) = w^{-1}(fg/\varphi)(\varphi^{-1}(w)) \in \dot{E}^r(\varphi(D))$. Since $\varphi(D)$ is bounded, we also have $w^{-1}(fg)(\varphi^{-1}(w)) \in \dot{E}^r(\varphi(D))$. 
Using Proposition \ref{EpDinvariantprop} again, we conclude that $fg/\varphi, fg \in \dot{E}^r(D)$. Part $(b)$ follows.
\proofend


\subsection{The spaces $\dot{L}^p(\Gamma)$}
Let $1 \leq p < \infty$. Let $\Gamma$ be a Carleson curve. We define $\dot{L}^p(\Gamma)$ as the set of all measurable functions on $\Gamma$ such that $|z - z_0|^{1 - \frac{2}{p}}h(z) \in L^p(\Gamma)$ for some (and hence all) $z_0 \in \C \setminus \Gamma$.
Note that 
\begin{align*}
& \dot{L}^p(\Gamma) \subset L^p(\Gamma), \qquad 2 \leq p < \infty,
	\\
& L^p(\Gamma) \subset \dot{L}^p(\Gamma), \qquad 1 \leq p \leq 2,
\end{align*}
If $h \in \dot{L}^p(\Gamma)$, then the value of the Cauchy integral $(\mathcal{C}h)(z)$ is well-defined for all $z \in \C \setminus \Gamma$. Indeed, if $1/p + 1/q = 1$, then the Carleson property of $\Gamma$ implies $\||\cdot - z|^{- \frac{2}{q}} \|_{L^q(\Gamma)} < \infty$; hence, by H\"older's inequality,
$$\int_\Gamma \frac{|h(z')|}{|z' - z|}|dz'|
\leq \big\||\cdot - z|^{1 - \frac{2}{p}} h \big\|_{L^p(\Gamma)} 
\big\||\cdot - z|^{- \frac{2}{q}}  \big\|_{L^q(\Gamma)} < \infty.$$
If $\Gamma$ is bounded, then $\dot{L}^p(\Gamma) = L^p(\Gamma)$. 

\begin{lemma}\label{Lpdotlemma}
Let $1 \leq p < \infty$ and let $\Gamma$ be a Carleson curve. Let $z_0 \in \C \setminus \Gamma$ and let $\varphi$ be given by (\ref{varphidef}).

\begin{enumerate}[$(a)$]
\item The map $\Phi$ defined for $h \in \dot{L}^p(\Gamma)$ by 
\begin{align}\label{Phidef}
(\Phi h)(w) = w^{-1} h(\varphi^{-1}(w)) 
\end{align}
is a bijection $\dot{L}^p(\Gamma) \to L^p(\varphi(\Gamma))$ and
$$\||\cdot - z_0|^{1-\frac{2}{p}} h\|_{L^p(\Gamma)} = \|\Phi h\|_{L^p(\varphi(\Gamma))}$$
for all $h \in \dot{L}^p(\Gamma)$.

\item If $h \in \dot{L}^p(\Gamma)$, then
\begin{align}\label{mathcalCinv}
(\mathcal{C}_\Gamma h)(z) = (\Psi^{-1} \mathcal{C}_{\varphi(\Gamma)} \Phi h)(z) \quad \text{for all $z \in \C \setminus \Gamma$},
\end{align}
where $\Psi$ acts on a function Ê$f:\hat{\C}\setminus \Gamma \to \C$ by $(\Psi f)(w) = w^{-1} f(\varphi^{-1}(w))$.
\end{enumerate}
\end{lemma}
\proofbegin
$(a)$ If $h \in \dot{L}^p(\Gamma)$, then the change of variables $w = \varphi(z)$ implies
\begin{align} \nonumber
\||z-z_0|^{1-\frac{2}{p}} h(z)\|_{L^p(\Gamma)}^p
& = \int_\Gamma |z-z_0|^{p-2} |h(z)|^p |dz|
	\\ \nonumber
& = \int_{\varphi(\Gamma)} |w^{-1} h(\varphi^{-1}(w))|^p |dw|
= \|\Phi h\|_{L^p(\varphi(\Gamma))}^p.
\end{align}
We infer that $\Phi$ is a bijection $\dot{L}^p(\Gamma) \to L^p(\varphi(\Gamma))$ with inverse given by $(\Phi^{-1}H)(z) = \varphi(z) H(\varphi(z))$. This proves $(a)$.

$(b)$  If $h \in \dot{L}^p(\Gamma)$, $z \in \C \setminus \Gamma$ and $w = \varphi(z)$, then
\begin{align}\nonumber
(\mathcal{C}_{\varphi(\Gamma)} \Phi h)(w)
& = \frac{1}{2\pi i} \int_{\varphi(\Gamma)} \frac{(h \circ \varphi^{-1})(w')}{w' - w}\frac{dw'}{w'}
	\\ 
& = \frac{z - z_0}{2\pi i} \int_{\Gamma} \frac{h(z')}{z' - z} dz'
= (z-z_0) (\mathcal{C}_\Gamma h)(z)
= \Psi(\mathcal{C}_\Gamma h)(w),
\end{align}
which proves $(b)$.
\proofend

For each $z_0 \in \C \setminus \Gamma$, we define a norm on $\dot{L}^p(\Gamma)$ by
\begin{align}\label{dotLpnorm}
\|h\|_{\dot{L}^p(\Gamma), z_0} = \||\cdot - z_0|^{1 - \frac{2}{p}} h\|_{L^p(\Gamma)}.
\end{align}
The space $\dot{L}^p(\Gamma)$ with the norm $\|\cdot\|_{\dot{L}^p(\Gamma), z_0}$ is nothing but the weighted space $L^p(\Gamma, w)$ with $w(z) =  |z - z_0|^{1 - \frac{2}{p}}$.
Different choices of $z_0 \in \C \setminus \Gamma$ induce different norms on $\dot{L}^p(\Gamma)$, but these norms are all equivalent. We say that an operator $T$ on $\dot{L}^p(\Gamma)$ is bounded if it is bounded with respect to one (and hence all) of these norms.

\begin{lemma}
Let $\Gamma \in \mathcal{J}$.
\begin{enumerate}[$(a)$]
\item $\dot{L}^p(\Gamma) \subset \dot{L}^r(\Gamma)$ whenever $1 \leq r \leq p < \infty$.
\item Suppose $p,q,r \in [1, \infty)$ satisfy $1/p + 1/q = 1/r$. If $f \in \dot{L}^p(\Gamma)$ and $g \in \dot{L}^q(\Gamma)$, then the functions $zf(z)g(z)$ and $fg$ belong to $\dot{L}^r(\Gamma)$.
\end{enumerate}
\end{lemma}
\proofbegin
$(a)$ The result is immediate from the definitions if $\infty \notin \Gamma$. Thus suppose $\infty \in \Gamma$. Let $z_0 \in \C \setminus \Gamma$ and define $\varphi$ as in (\ref{varphidef}). If $h \in \dot{L}^p(\Gamma)$, then Lemma \ref{Lpdotlemma} implies that $w^{-1}h(\varphi^{-1}(w)) \in \dot{L}^p(\varphi(\Gamma))$; since $\dot{L}^p(\varphi(\Gamma)) \subset \dot{L}^r(\varphi(\Gamma))$, another application of Lemma \ref{Lpdotlemma} shows thatÊ $h \in \dot{L}^r(\Gamma)$. This proves $(a)$. 

$(b)$ Suppose $p,q,r \in [1, \infty)$ satisfy $1/p + 1/q = 1/r$. Let $f \in \dot{L}^p(\Gamma)$ and $g \in \dot{L}^q(\Gamma)$. 
Suppose $\infty \notin \Gamma$. Then, by H\"older's inequality,
$$\|fg\|_{L^r(\Gamma)} \leq \|f\|_{L^p(\Gamma)} \|g\|_{L^q(\Gamma)} < \infty.$$
Hence $fg \in L^r(\Gamma) = \dot{L}^r(\Gamma)$. This gives $(b)$ when $\infty \notin \Gamma$.

If $\infty \in \Gamma$, then pick $z_0 \in \C \setminus \Gamma$ and define $\varphi$ as in (\ref{varphidef}). By Lemma \ref{Lpdotlemma}, $w^{-1} f(\varphi^{-1}(w)) \in L^p(\varphi(\Gamma))$ and $w^{-1} g(\varphi^{-1}(w)) \in L^q(\varphi(\Gamma))$. Hence, by the preceding paragraph, 
$w^{-2}(fg)(\varphi^{-1}(w)) = w^{-1}(fg/\varphi)(\varphi^{-1}(w)) \in L^r(\varphi(\Gamma))$. Since $\varphi(\Gamma)$ is bounded, we also have $w^{-1}(fg)(\varphi^{-1}(w)) \in L^r(\varphi(\Gamma))$. 
Using Lemma \ref{Lpdotlemma} again, we conclude that $fg/\varphi, fg \in \dot{L}^r(\Gamma)$. Part $(b)$ follows.
\proofend

\subsection{The Cauchy singular operator}
Theorem \ref{Sth} can be used to establish boundedness of the Cauchy singular operator $\mathcal{S}_\Gamma$ on $\dot{L}^p(\Gamma)$ if $1 < p < \infty$ and $\Gamma$ is Carleson.

\begin{proposition}\label{SLpdotprop}
Let $1 < p < \infty$ and let $\Gamma$ be a Carleson curve. 
Then $\mathcal{S}_\Gamma$ generates a bounded operator $\tilde{\mathcal{S}}_\Gamma$ on $\dot{L}^p(\Gamma)$. Moreover, if $h \in \dot{L}^p(\Gamma)$, then the limit in (\ref{cauchysingulardef}) exists and $(\mathcal{S}_\Gamma h)(z) = (\tilde{\mathcal{S}}_\Gamma h)(z)$ for a.e. $z \in \Gamma$.
\end{proposition}
\proofbegin
Let $z_0 \notin \Gamma$ and let $w(z) = |z - z_0|^{1 - \frac{2}{p}}$. The result follows from Theorem \ref{Sth} if we can show that $w \in A_p(\Gamma)$. If $p =2$, this is an immediate consequence of $\Gamma$ being Carleson. Thus suppose $p \neq 2$.

Define $q$ by $1/p + 1/q = 1$ and let $z \in \Gamma$. If $0 < r \leq \frac{|z-z_0|}{2}$ and $z' \in D(z,r)$, then
$$ \frac{|z - z_0|}{2}  \leq |z' - z_0| \leq \frac{3|z - z_0|}{2} .$$
Using the Carleson condition (\ref{carlesondef}) on the disk $D(z, r)$, we find that there exists a constant $C_\alpha > 0$ depending only on $\alpha$ such that
\begin{subequations}\label{smallrest}
\begin{align}\label{smallrest1}
& \frac{1}{r} \int_{\Gamma \cap D(z, r)} |z' - z_0|^\alpha |dz'|
\leq
\frac{1}{r} \frac{3^\alpha |z-z_0|^\alpha}{2^\alpha} |\Gamma \cap D(z, r)|
\leq C_\alpha |z-z_0|^\alpha,
	\\ \label{smallrest2}
& \frac{1}{r} \int_{\Gamma \cap D(z, r)} |z' - z_0|^{-\alpha} |dz'|
\leq
\frac{1}{r} \frac{|z-z_0|^{-\alpha}}{2^{-\alpha}} |\Gamma \cap D(z, r)|
\leq C_\alpha |z-z_0|^{-\alpha},
\end{align}
\end{subequations}
whenever $0 < r \leq \frac{|z-z_0|}{2}$ and $\alpha > 0$.
This yields
\begin{align}\nonumber
& \sup_{z \in \Gamma} \sup_{0 < r \leq \frac{|z-z_0|}{2}} \bigg(\frac{1}{r} \int_{\Gamma \cap D(z, r)} w(z')^p |dz'|\bigg)^{1/p}\bigg(\frac{1}{r} \int_{\Gamma \cap D(z, r)} w(z')^{-q} |dz'|\bigg)^{1/q} 
	\\ \nonumber
= & \sup_{z \in \Gamma} \sup_{0 < r \leq \frac{|z-z_0|}{2}} \bigg(\frac{1}{r} \int_{\Gamma \cap D(z, r)} |z' - z_0|^{p -2} |dz'|\bigg)^{1/p}\bigg(\frac{1}{r} \int_{\Gamma \cap D(z, r)} |z' - z_0|^{\frac{q}{p} - 1} |dz'|\bigg)^{1/q} 
	\\ \label{muckpart1}
& \leq
C_p \sup_{z \in \Gamma} \bigg(|z-z_0|^{p - 2} \bigg)^{1/p}\bigg(|z-z_0|^{\frac{q}{p} - 1} \bigg)^{1/q} 
= C_p < \infty,
\end{align}
with $C_p > 0$ depending only on $p$.

On the other hand, if $R = |z - z_0| + r$, then the Carleson condition on the disk $D(z_0, 2^{1-n}R)$ yields
\begin{align}\nonumber
\frac{1}{r} \int_{\Gamma \cap D(z, r)} & |z' - z_0|^\alpha |dz'|
 =
\frac{1}{r} \sum_{n=1}^\infty \int_{\Gamma \cap D(z, r) \cap (D(z_0, 2^{1-n}R) \setminus D(z_0, 2^{-n}R))} |z' - z_0|^\alpha |dz'|
	\\\nonumber
& \leq
\frac{1}{r} \sum_{n=1}^\infty 2^{(1-n)\alpha }R^\alpha  |\Gamma \cap (D(z_0, 2^{1-n}R) \setminus D(z_0, 2^{-n}R))|
	\\ \label{bigrest1}
& \leq
\frac{C R^{1+ \alpha} }{r} \sum_{n=1}^\infty 2^{(1-n)\alpha }  2^{1-n}
\leq \frac{2C (|z - z_0| + r)^{1 + \alpha}}{r}
\leq C_\alpha r^\alpha
\end{align}
whenever $r > \frac{|z-z_0|}{2}$ and $\alpha > 0$. Similarly,
\begin{align}\nonumber
\frac{1}{r} \int_{\Gamma \cap D(z, r)}& |z' - z_0|^{-\alpha} |dz'|
 =
\frac{1}{r} \sum_{n=1}^\infty \int_{\Gamma \cap D(z, r) \cap (D(z_0, 2^{1-n}R) \setminus D(z_0, 2^{-n}R))} |z' - z_0|^{-\alpha} |dz'|
	\\ \nonumber
& \leq
\frac{1}{r} \sum_{n=1}^\infty 2^{n\alpha }R^{-\alpha}  |\Gamma \cap (D(z_0, 2^{1-n}R) \setminus D(z_0, 2^{-n}R))|
	\\Ê\label{bigrest2}
& \leq
\frac{C R^{1- \alpha} }{r} \sum_{n=1}^\infty 2^{n\alpha }  2^{1-n}
\leq \frac{C_\alpha (|z - z_0| + r)^{1 - \alpha}}{r}
\leq C_\alpha r^{- \alpha},
\end{align}
whenever $r > \frac{|z-z_0|}{2}$ and $0 < \alpha < 1$.
If $p > 2$, we apply (\ref{bigrest1}) with $\alpha = p - 2$ and (\ref{bigrest2}) with $\alpha = 1- \frac{q}{p}$. If $1 < p < 2$, we apply (\ref{bigrest2}) with $\alpha = 2- p$ and (\ref{bigrest1}) with $\alpha = \frac{q}{p} - 1$. This yields
\begin{align}\nonumber
& \sup_{z \in \Gamma} \sup_{r > \frac{|z-z_0|}{2}} \bigg(\frac{1}{r} \int_{\Gamma \cap D(z, r)} w(z')^p |dz'|\bigg)^{1/p}\bigg(\frac{1}{r} \int_{\Gamma \cap D(z, r)} w(z')^{-q} |dz'|\bigg)^{1/q} 
	\\ \label{muckpart2}
& \leq C \sup_{z \in \Gamma} \sup_{r > \frac{|z-z_0|}{2}} \big(r^{p -2} \big)^{1/p}
\big(r^{\frac{q}{p} - 1} \big)^{1/q} 
 = C < \infty.
\end{align}
It follows from (\ref{muckpart1}) and (\ref{muckpart2}) that $w(z)$ satisfies the Muckenhoupt condition (\ref{muckenhouptcondition}).
\proofend

Our next objective is to determine how $\mathcal{S}_\Gamma$ transforms under the change of variables Ê$w = 1/(z - z_0)$. We need the following lemma.

\begin{lemma}\label{annuluslemma}
  Let $\Gamma$ be a Carleson jump contour. 
  Let $g(\epsilon)$ be a nondecreasing continuous function of Ê$\epsilon \geq 0$ such that $g(0) = 0$.
If $z \in \Gamma \cap \C$ is a point at which $\Gamma$ has a two-sided tangent, then the following limit exists and equals zero:
\begin{align}\label{limepsilong}
\lim_{\epsilon \to 0} \int_{\Gamma \cap [D(z, \epsilon(1 + g(\epsilon))) \setminus D(z, \epsilon(1 - g(\epsilon)))]}\frac{|dz'|}{|z' - z|} = 0.
\end{align}
\end{lemma}
\proofbegin
Without loss of generality, we assume that $\Gamma \in \mathcal{J}$.
Let $z \in \Gamma \cap \C$ Êbe a point at which  $\Gamma$ has a two-sided tangent. Let $\gamma(s)$, $-s_0 < s < s_0$, be an arclength parametrization of $\Gamma$ in a neighborhood of $\gamma(0) = z$. Suppose without loss of generality that $\gamma'(0) = 1$. Then 
$$\gamma(s) = z + s + o(|s|), \qquad s \to 0.$$
For each $r \in (0,1/2]$, choose $\delta(r) \in (0, s_0)$ such that
\begin{align}\label{oestimate}
|o(|s|)| < r |s| \quad \text{for}Ê\quad |s| \leq \delta(r).
\end{align}
We may assume that $\delta(r)$ is a nondecreasing function of $r > 0$. Replacing $\delta(r)$ with $\int_0^r \delta(t)dt \leq \delta(r)$ if necessary (note that all nondecreasing functions are measurable), we may assume that $\delta(r)$ is a continuous strictly increasing function such that $ \lim_{r \to 0^+} \delta(r) = 0$.
Let $\gamma_{1/2}$ denote the subarc
$$\gamma_{1/2} = \{\gamma(s) | |s| \leq \delta(1/2)\}$$
and let $a, b$ be the endpoints of $\gamma_{1/2}$. The set $(\Gamma \setminus \gamma_{1/2}) \cup \{a, b\}$ is compact. Let $\mu$ be the minimum of the continuous function $|\cdot - z|$ on this set. Then $\Gamma \cap D(z, \mu) \subset \gamma_{1/2}$. 
Fix $r\in (0,1/2]$ such that $\delta(r) < \mu$. We claim that 
\begin{align}\label{zznn}
|s| < \frac{\delta(r)}{2(1-r)},
\end{align}
whenever $\gamma(s) \in \gamma_{1/2}$ and $|\gamma(s) - z| \leq  \frac{\delta(r)}{2}$.
Indeed, suppose $|s| \leq \delta(1/2)$ is such that $|\gamma(s) - z| \leq \frac{\delta(r)}{2}$. Then (\ref{oestimate}) implies
$$\frac{|s|}{2}Ê< |s + o(|s|)| = |\gamma(s) - z| \leq \frac{\delta(r)}{2}.$$
Thus $|s| < \delta(r)$, so another application of (\ref{oestimate}) yields
$$(1 - r) |s| < |s + o(|s|)| = |\gamma(s) - z| \leq \frac{\delta(r)}{2},$$
which proves (\ref{zznn}).
On the other hand, since $s$ is an arclength parameter,
\begin{align}\label{zz} 
  |s| \geq \frac{\delta(r)}{2}\bigg(1 - g\bigg(\frac{\delta(r)}{2}\bigg)\bigg),
\end{align}
whenever $|\gamma(s) - z| \geq \frac{\delta(r)}{2}(1 - g(\frac{\delta(r)}{2}))$.

For $\epsilon > 0$, we define the closed annulus $K(z, \epsilon)$ by
$$K(z, \epsilon) = \overline{D(z, \epsilon)} \setminus D(z, \epsilon(1 - g(\epsilon))).$$
Then the set $\{s \in [0,\delta(1/2)] | \gamma(s) \in K(z, \delta(r)/2)\}$ is closed. Let $s_+ \geq 0$Ê and $s_- \geq 0$ denote the largest and smallest elements of this set, respectively. Clearly,
$$|\gamma(s_+) - z| \leq \frac{\delta(r)}{2}, \qquad
|\gamma(s_-) - z| \geq \frac{\delta(r)}{2}\bigg(1 - g\bigg(\frac{\delta(r)}{2}\bigg)\bigg).$$
Hence, by (\ref{zznn}) and (\ref{zz}),
$$\frac{\delta(r)}{2}\bigg(1 - g\bigg(\frac{\delta(r)}{2}\bigg)\bigg) \leq s_- \leq s_+ < \frac{\delta(r)}{2(1-r)}.$$
Thus
$$\bigg|\gamma([0,\delta(1/2)]) \cap K\bigg(z, \frac{\delta(r)}{2}\bigg)\bigg| \leq s_+ - s_-
< \frac{\delta(r)}{2} \bigg( \frac{r}{1-r} + g\bigg(\frac{\delta(r)}{2}\bigg)\bigg).$$
A similar argument shows that
$$\bigg|\gamma([-\delta(1/2), 0]) \cap K\bigg(z, \frac{\delta(r)}{2}\bigg)\bigg| < \frac{\delta(r)}{2} \bigg( \frac{r}{1-r} + g\bigg(\frac{\delta(r)}{2}\bigg)\bigg).$$
Consequently, for all small enough $r > 0$,
\begin{align}\label{GammaK}
\bigg|\Gamma \cap K\bigg(z, \frac{\delta(r)}{2}\bigg)\bigg|
= \bigg|\gamma_{1/2} \cap K\bigg(z, \frac{\delta(r)}{2}\bigg)\bigg|
< F(r)
\end{align}
where the function 
$$F(r) = \delta(r) \bigg( \frac{r}{1-r} + g\bigg(\frac{\delta(r)}{2}\bigg)\bigg)$$
satisfies
$$\lim_{r \to 0^+} \frac{F(r)}{\delta(r)} = 0.$$
Given $\epsilon>0$ small enough, there exists a unique $r = r(\epsilon) > 0$ such that $\delta(r)/2 = \epsilon$. 
It follows that
\begin{align*}
\int_{\Gamma \cap K(z,\epsilon)}\frac{|dz'|}{|z' - z|} 
\leq \frac{|\Gamma \cap K(z,\epsilon)|}{\epsilon(1 - g(\epsilon))}
< \frac{2 F(r)}{\delta(r) (1 - g(\epsilon))}
\to 0
\end{align*}
as $\epsilon \to 0$. This proves that
\begin{align}\label{limepsilong2}
\lim_{\epsilon \to 0} \int_{\Gamma \cap [D(z, \epsilon) \setminus D(z, \epsilon(1 - g(\epsilon)))]}\frac{|dz'|}{|z' - z|} = 0.
\end{align}
Equation (\ref{limepsilong}) follows from (\ref{limepsilong2}) by changing variables $\tilde{\epsilon} = \epsilon(1 + g(\epsilon))$ in the left-hand side of (\ref{limepsilong}) and noting that $\epsilon(1 - g(\epsilon)) = \tilde{\epsilon}(1 - \tilde{g}(\tilde{\epsilon}))$, where $\tilde{g}(\tilde{\epsilon}) = \frac{2g(\epsilon)}{1 + g(\epsilon)}$ is a continuous nondecreasing function of $\tilde{\epsilon}$.
\proofend

\begin{proposition}
Let $1 < p < \infty$ and let $\Gamma$ be a Carleson jump contour. Let $z_0 \in \C \setminus \Gamma$ and let $\varphi(z) = 1/(z-z_0)$. Let $\Phi:\dot{L}^p(\Gamma) \to L^p(\varphi(\Gamma))$ be the bijection defined in (\ref{Phidef}). Then
\begin{align}\label{SGammainvariance}
\mathcal{S}_\Gamma h = \Phi^{-1} \mathcal{S}_{\varphi(\Gamma)} \Phi h \quad \text{a.e. on $\Gamma$}
\end{align}
for every $h \in \dot{L}^p(\Gamma)$. In other words, the following diagram commutes:
$$
\xymatrix{\ar[d]_{\mathcal{S}_\Gamma} \dot{L}^p(\Gamma) \ar[r]^{\Phi} & L^p(\varphi(\Gamma)) \ar[d]^{\mathcal{S}_{\varphi(\Gamma)}} \\
 \dot{L}^p(\Gamma) \ar[r]^{\Phi} &  L^p(\varphi(\Gamma))}
$$
\end{proposition}
\proofbegin
We will show that $\mathcal{S}_\Gamma h = \Phi^{-1} \mathcal{S}_{\varphi(\Gamma)} \Phi h$ a.e. on $\Gamma$ whenever $h \in C_0^\infty(\Gamma)$. Since $C_0^\infty(\Gamma)$ is dense in $\dot{L}^p(\Gamma)$ and the operators $\mathcal{S}_\Gamma$ and $\Phi^{-1} \mathcal{S}_{\varphi(\Gamma)} \Phi$ are bounded on $\dot{L}^p(\Gamma)$ by Lemma \ref{Lpdotlemma} and Proposition \ref{SLpdotprop}, this will prove (\ref{SGammainvariance}).

Let $h \in C_0^\infty(\Gamma)$. A change of variables shows that
\begin{align}\label{SGammaints}
\frac{1}{\pi i} \int_{\Gamma \setminus D(z, \epsilon)} \frac{h(z')}{z' - z}  dz'
= \frac{\varphi(z)}{\pi i} \int_{\varphi(\Gamma) \setminus \varphi(D(z, \epsilon))} \frac{(\Phi h)(w')}{w' - \varphi(z)}  dw'
\end{align}
for all $z \in \Gamma \cap \C$ and $\epsilon > 0$. As $\epsilon \to 0$, the left-hand side of (\ref{SGammaints}) tends to $(\mathcal{S}_\Gamma h)(z)$ for a.e. $z \in \Gamma$. It remains to prove that the right-hand side of (\ref{SGammaints}) tends to 
\begin{align}\label{PhiinvS}
(\Phi^{-1} \mathcal{S}_{\varphi(\Gamma)} \Phi h)(z)
= \varphi(z) (\mathcal{S}_{\varphi(\Gamma)} \Phi h)(\varphi(z))
\end{align}
for a.e. $z \in \Gamma$ as $\epsilon \to 0$. The proof of this fact is complicated by the fact that,  in general, the disk $\varphi(D(z, \epsilon))$ is not centered at $\varphi(z)$.

Let $z \in \Gamma$ and let $0 < \epsilon < |z-z_0|$. Then
$$\varphi(D(z,\epsilon)) = D\bigg( \frac{\bar{z} - \bar{z}_0}{|z - z_0|^2 - \epsilon^2}, \tilde{\epsilon}\bigg)$$
and
\begin{align}\label{nesteddisks}
D\bigg(\varphi(z), \tilde{\epsilon}\bigg(1 - \frac{\epsilon}{|z - z_0|}\bigg)\bigg) \subset \varphi(D(z,\epsilon))  \subset D\bigg(\varphi(z), \tilde{\epsilon}\bigg(1 + \frac{\epsilon}{|z - z_0|}\bigg)\bigg),
\end{align}
where
$$\tilde{\epsilon} = \frac{\epsilon}{|z - z_0|^2 - \epsilon^2}.$$
Noting that
$$1 \pm \frac{\epsilon}{|z - z_0|} = 1 \pm g(\tilde{\epsilon})$$
where 
$$g(\tilde{\epsilon}) = \frac{\sqrt{1 + 4|z - z_0|^2 \tilde{\epsilon}^2} - 1}{2|z - z_0|\tilde{\epsilon}},$$ 
we can write (\ref{nesteddisks}) as 
\begin{align}\label{nesteddisks2}
D(\varphi(z), \tilde{\epsilon}(1 - g(\tilde{\epsilon}))) \subset \varphi(D(z,\epsilon))  \subset D(\varphi(z), \tilde{\epsilon}(1 + g(\tilde{\epsilon}))).
\end{align}

The function $(\Phi h)(w) = w^{-1}h(w^{-1} + z_0)$ is the restriction to $\varphi(\Gamma)$Ê of a smooth function which approaches zero as $w \to \infty$. Hence there exists an $M >0$ such that $|(\Phi h)(w)| \leq M$ for all $w \in \varphi(\Gamma)$. We estimate
\begin{align}\nonumber
\bigg| \int_{\varphi(\Gamma) \setminus \varphi(D(z, \epsilon))} &\frac{(\Phi h)(w')}{w' - \varphi(z)}  dw' - \int_{\varphi(\Gamma) \setminus D(\varphi(z),  \tilde{\epsilon}(1 + g(\tilde{\epsilon})))} \frac{(\Phi h)(w')}{w' - \varphi(z)}  dw'\bigg|
	\\ \nonumber
& = 
\bigg| \int_{\varphi(\Gamma) \cap [D(\varphi(z),  \tilde{\epsilon}(1 + g(\tilde{\epsilon}))) \setminus \varphi(D(z, \epsilon))]} \frac{(\Phi h)(w')}{w' - \varphi(z)}  dw' \bigg|
	\\ \label{intdifferenceestimate}
& \leq M 
 \int_{\varphi(\Gamma) \cap [D(\varphi(z),  \tilde{\epsilon}(1 + g(\tilde{\epsilon}))) \setminus D(\varphi(z),  \tilde{\epsilon}(1 - g(\tilde{\epsilon})))]} \frac{|dw'|}{|w' - \varphi(z)|}.
 \end{align}
Being locally rectifiable, the Carleson jump contour $\varphi(\Gamma)$ has a  two-sided tangent at almost every point. Hence, by Lemma \ref{annuluslemma}, the limit of the right-hand side of (\ref{intdifferenceestimate}) as $\epsilon \to 0$ exists and equals zero for a.e. $z \in \Gamma$.
On the other hand, by Proposition \ref{SLpdotprop}, the limit
$$\lim_{\epsilon \to 0}Ê\frac{1}{\pi i} \int_{\varphi(\Gamma) \setminus D(\varphi(z), \tilde{\epsilon}(1 + g(\tilde{\epsilon})))} \frac{(\Phi h)(w')}{w' - \varphi(z)}  dw'$$
exists and equals $(S_{\varphi(\Gamma)}\Phi h)(\varphi(z))$ for a.e. $z \in \Gamma$. It follows that
$$\lim_{\epsilon \to 0} \frac{\varphi(z)}{\pi i} \int_{\varphi(\Gamma) \setminus \varphi(D(z, \epsilon))} \frac{(\Phi h)(w')}{w' - \varphi(z)}  dw' = \varphi(z) (S_{\varphi(\Gamma)}\Phi h)(\varphi(z))$$
for a.e. $z \in \Gamma$. This completes the proof.
\proofend

\section{Cauchy integrals over Carleson jump contours}\nequation\label{cauchysec}
The following two theorems generalize Theorems \ref{Jordanth1} and \ref{Jordanth2} to the case where $\Gamma$ is a general Carleson jump contour.

\begin{theorem}\label{jumpth1}
Let $\Gamma \subset \hat{\C}$ be a Carleson jump contour and let $D_\pm \subset \hat{\C}$ be the associated open sets such that $\partial D_+ = -\partial D_- = \Gamma$.
Let $\mathcal{C}$ denote the Cauchy integral operator defined in (\ref{cauchyintegraldef}).
\begin{enumerate}[$(a)$]
\item Let $1 \leq p < \infty$. If $f \in \dot{E}^p(D_+)$, then the nontangential limits of $f(z)$ as $z$ approaches the boundary exist a.e. on $\Gamma$; if $f_+(z)$ denotes the boundary function, then $f_+ \in \dot{L}^p(\Gamma)$ and
\begin{align}\label{jumpCfplus}
(\mathcal{C} f_+)(z) = \begin{cases} f(z), & z \in D_+, \\
0, & z \in D_-.
\end{cases}
\end{align}
If $f \in \dot{E}^p(D_-)$, then the nontangential limits of $f(z)$ as $z$ approaches the boundary exist a.e. on $\Gamma$. If $f_-(z)$ denotes the boundary function, then $f_- \in \dot{L}^p(\Gamma)$ and
\begin{align}\label{jumpCfminus}
(\mathcal{C} f_-)(z) = \begin{cases} 0, & z \in D_+, \\
 - f(z), & z \in D_-.
\end{cases}
\end{align}
In particular, $f = \mathcal{C}(f_+ - f_-)$ for all $f \in \dot{E}^p(D_+ \cup D_-)$.

\item Let $1 < p < \infty$. Then the Cauchy singular operator $\mathcal{S}_\Gamma:\dot{L}^p(\Gamma) \to \dot{L}^p(\Gamma)$ defined in (\ref{cauchysingulardef}) satisfies $\mathcal{S}_\Gamma^2 = I$. Moreover,
if $h \in \dot{L}^p(\Gamma)$, then 
\begin{align}\label{jumpChD}
\mathcal{C}h|_{D_+} \in \dot{E}^p(D_+), \qquad \mathcal{C}h|_{D_-} \in \dot{E}^p(D_-).
\end{align}
\end{enumerate}
\end{theorem}

Theorem \ref{jumpth1} implies that if $\Gamma$ is a Carleson jump contour and $h \in \dot{L}^p(\Gamma)$ for some $1 < p < \infty$, then the left and right nontangential boundary values of $\mathcal{C}h$, which we denote by $\mathcal{C}_+ h$ and $\mathcal{C}_- h$, lie in $\dot{L}^p(\Gamma)$. This allows us to define two linear operators $\mathcal{C}_\pm:h \mapsto \mathcal{C}_\pm h$ on $\dot{L}^p(\Gamma)$.

\begin{theorem}\label{jumpth2}
Let $1 < p < \infty$ and let $\Gamma \subset \hat{\C}$ be a Carleson jump contour.
Then $\mathcal{C}_\pm$ are bounded operators on $\dot{L}^p(\Gamma)$ with the following properties:

\begin{itemize}
\item The Sokhotski-Plemelj formulas 
\begin{align}\label{plemelj}
\mathcal{C}_+ = \frac{1}{2}(I + \mathcal{S}_\Gamma), \qquad \mathcal{C}_- = \frac{1}{2}(-I + \mathcal{S}_\Gamma),
\end{align}
are valid. 

\item $\mathcal{C}_\pm$ are orthogonal projections on $\dot{L}^p(\Gamma)$ in the sense that
$$\dot{L}^p(\Gamma) = \mathcal{C}_+\dot{L}^p(\Gamma) \oplus \mathcal{C}_-\dot{L}^p(\Gamma)$$ 
and
$$\mathcal{C}_+ - \mathcal{C}_- = \I, \qquad \mathcal{C}_+^2 = \mathcal{C}_+, \qquad \mathcal{C}_-^2 = -\mathcal{C}_-, \qquad \mathcal{C}_+\mathcal{C}_- = \mathcal{C}_-\mathcal{C}_+ = 0.$$

\item If $h = \mathcal{C}_+h - \mathcal{C}_-h \in \dot{L}^p(\Gamma)$, then
\begin{align}\label{ChDChD}
(\mathcal{C}h)|_{D_+} = (\mathcal{C}\mathcal{C}_+h)|_{D_+} \in \dot{E}^p(D_+), \qquad
(\mathcal{C}h)|_{D_-} = -(\mathcal{C}\mathcal{C}_-h)|_{D_-} \in \dot{E}^p(D_-).
\end{align}

\item The map $h \mapsto (\mathcal{C}h)|_{D_+}$ is a bijection  $\mathcal{C}_+\dot{L}^p(\Gamma) \to \dot{E}^p(D_+)$ with inverse $f \mapsto f_+$. 

\item The map $h \mapsto (\mathcal{C}h)|_{D_-}$ is a bijection  $\mathcal{C}_-\dot{L}^p(\Gamma) \to \dot{E}^p(D_-)$ with inverse $f \mapsto -f_-$.
\end{itemize}
\end{theorem}

In the special case of a jump contour $\Gamma$ consisting of a single rectifiable Jordan curve, Theorems \ref{jumpth1} and \ref{jumpth2} reduce to Theorems \ref{Jordanth1} and \ref{Jordanth2}, respectively.

\subsection{Proof of Theorem \ref{jumpth1}}
\subsubsection{Proof of $(a)$}
Suppose first that $\infty \notin \Gamma$, so that $\Gamma \subset \C$ is bounded. 
Let $f \in \dot{E}^p(D_+)$. 
Represent $\Gamma$ as the union of finitely many arcs each pair of which have at most endpoints in common. If $z \in \Gamma$ is not one of these finitely many endpoints, then $z$ belongs to $\partial D_j^+$ for exactly one component $D_j^+$ of $D_+$. Since Theorem \ref{Jordanth1} implies that $f|_{D_j^+}$ has nontangential limits a.e. on $\partial D_j^+$, it follows that $f$ has nontangential limits a.e. on $\Gamma$.
Another application of Theorem \ref{Jordanth1} shows that $f_+|_{\partial D_j^+} \in L^p(\partial D_j^+)$ for each $j$. Hence 
$f_+ \in L^p(\Gamma) = \dot{L}^p(\Gamma)$. 

Now suppose $z \in D_{k}^+$ for some $1 \leq k \leq n$. Since $z$ lies in the region exterior to $\partial D_j^+$ for each $j \neq k$, Theorem \ref{Jordanth1} yields
$$(\mathcal{C}f_+)(z) = \frac{1}{2\pi i} \int_{\partial D_k^+} \frac{f_+(z')}{z' - z}dz'
 + \frac{1}{2\pi i} \sum_{j \neq k} \int_{\partial D_j^+} \frac{f_+(z')}{z' - z}dz'
 = f(z).$$
If $z \in D_-$, then $z$ lies in the region exterior to $\partial D_j^+$ for every $j$, so a similar computation implies $(\mathcal{C}f_+)(z) = 0$. This proves (\ref{jumpCfplus}).
Similar arguments apply when $f \in \dot{E}^p(D_-)$. This proves $(a)$ in the case when $\Gamma$ is bounded. 

Suppose now that $\infty \in \Gamma$. Pick $z_0 \in D_-$ and let $\varphi(z) = 1/(z-z_0)$. Let $f \in \dot{E}^p(D_+)$.
Then $\infty \notin \varphi(\Gamma)$ and $\varphi(\Gamma)$ is a Carleson jump contour by Proposition \ref{jumpprop}.
Let $F(w) = w^{-1} f(\varphi^{-1}(w))$. Then $F \in \dot{E}^p(\varphi(D_+))$ and $f \circ \varphi^{-1} \in E^p(\varphi(D_+))$ by Proposition \ref{EpDinvariantprop}. Since $\varphi(\Gamma)$ is bounded, the result of the preceding paragraph implies that the nontangential boundary values of $f \circ \varphi^{-1}$ exist a.e. on $\varphi(\Gamma)$. It follows that the nontangential boundary values of $f$ exist a.e. on $\Gamma$ and $(f \circ \varphi^{-1})_+ = f_+ \circÊ\varphi^{-1}$ a.e. on $\varphi(\Gamma)$.
Furthermore, since $F \in \dot{E}^p(\varphi(D_+))$, we have $F_+ \in \dot{L}^p(\varphi(\Gamma))$, which by Lemma \ref{Lpdotlemma} implies that $f_+ \in \dot{L}^p(\Gamma)$.
We also have
\begin{align}\label{CFplus}
(\mathcal{C}_{\varphi(\Gamma)} F_+)(w) = \begin{cases} F(w), & w \in \varphi(D_+), \\
0, & w \in \varphi(D_-),
\end{cases}
\end{align}
which in view of Lemma \ref{Lpdotlemma} yields (\ref{jumpCfplus}).
Similar arguments apply when $f \in \dot{E}^p(D_-)$. This proves $(a)$.

\subsubsection{A convergence lemma}
For the proof of $(b)$, we need the following lemma.

\begin{lemma}\label{Epconvergencelemma}
Let $1 \leq p < \infty$. Let $\Gamma \subset \C$ be a rectifiable Jordan curve and let $D_+$ and $D_-$ be the interior and exterior components of $\hat{\C} \setminus \Gamma$. Suppose $h \in L^p(\Gamma)$.

\begin{enumerate}[$(i)$]
\item If $\{f_n\}_1^\infty$ is a sequence of functions in $E^p(D_+)$ such that $f_{n+} \to h$ in $L^p(\Gamma)$, then there exists a function $f \in E^p(D_+)$ such that $f_n \to f$ uniformly on compact subsets of $D_+$ and $f_+ = h$.  

\item If $\{f_n\}_1^\infty$ is a sequence of functions in $E^p(D_-)$ such that $f_{n-} \to h$ in $L^p(\Gamma)$, then there exists a function $f \in E^p(D_-)$ such that $f_n \to f$ uniformly on compact subsets of $D_-$ and $f_- = h$. If $\{f_n\}_1^\infty \subset \dot{E}^p(D_-)$, then $f \in \dot{E}^p(D_-)$.
\end{enumerate}
\end{lemma}
\proofbegin
Part $(i)$ is a consequence of Theorem 17.2 in Chapter III of \cite{P1956}. 
In order to prove $(ii)$, let $\{f_n\}_1^\infty$ be a sequence of functions in $E^p(D_-)$ such that $f_{n-} \to h$ in $L^p(\Gamma)$. Let $z_0\in D_+$ and let $\varphi(z) = 1/(z-z_0)$. Then $h \circ \varphi^{-1} \in L^p(\varphi(\Gamma))$ and Proposition \ref{EpDinvariantprop} implies that $f_n \circ \varphi^{-1} \in E^p(\varphi(D_-))$ for each $n$. Assuming for simplicity that both $\Gamma$Ê and $\varphi(\Gamma)$ are oriented counterclockwise, we have $(f_n \circ \varphi^{-1})_+ = f_{n-} \circ \varphi^{-1}$, and so
\begin{align*}
\|(f_n \circ \varphi^{-1})_+ - h \circ \varphi^{-1}\|_{L^p(\varphi(\Gamma))}^p
= \int_{\varphi(\Gamma)} |f_{n-}(\varphi^{-1}(w)) - h(\varphi^{-1}(w))|^p |dw|
	\\
= \int_\Gamma |f_{n-}(z) - h(z)|^p \frac{|dz|}{|z-z_0|^2}
\leq C \|f_{n-} - h\|_{L^p(\Gamma)}^p \to 0 \quad \text{as} \quad n \to \infty.
\end{align*}
Hence, by $(i)$, there exists a function $g \in E^p(\varphi(D_-))$ such that
$f_n \circ \varphi^{-1} \to g$ uniformly on compact subsets of $\varphi(D_-)$ and $g_+ = h \circ \varphi^{-1}$. Letting $f = g \circ \varphi$, we infer that $f \in E^p(D_-)$,  that $f_n \to f$ uniformly on compact subsets of $D_-$, and that $f_- = h$. If $\{f_n\}_1^\infty \subset \dot{E}^p(D_-)$, then each $f_n$ vanishes at $\infty$. Hence $f$ vanishes at $\infty$, and so $f \in \dot{E}^p(D_-)$.
\proofend

\subsubsection{Proof of $(b)$}
 Suppose $1 < p < \infty$ and $h \in \dot{L}^p(\Gamma)$. 
We first assume that $\infty \notin \Gamma$. Switching the orientation of $\Gamma$ if necessary, we may suppose that $\infty \in D_-$. 
Let $R(\Gamma)$ be the set of all rational functions with no poles on $\Gamma$.
Every function $r \in R(\Gamma)$ can be written as $r = r^+ + r^-$, where $r^+$ is analytic in $D_+$, $r^-$ is analytic in $D_-$, and $r^-$ vanishes at infinity.
That is, $r^+ \in \dot{E}^p(D_+)$ and $r^- \in \dot{E}^p(D_-)$.
We claim that
\begin{align}\label{SGammar}
  \mathcal{S}_\Gamma r^+ = r^+, \qquad \mathcal{S}_\Gamma r^- = -r^-.
\end{align}
Indeed, if $\Gamma$ consists of a single Carleson Jordan curve, then (\ref{SGammar}) follows from Lemma 6.5 of \cite{BK1997}. If $\Gamma$ is the union of multiple Carleson Jordan curves $\{\partial D_j^+\}_1^n$, then we write $r^- = \sum_{j=1}^n r_j^-$ where $r_j^-$ is analytic outside $D_j^+$ and $r_j^-(\infty) = 0$ for each $j$. Let $\chi_j$ be the characteristic function of $\partial D_j^+$.
Decomposing $r^+$ and $r_j^-$ into partial fractions and using that (\ref{SGammar}) is valid in the case when $\Gamma$ is a Carleson Jordan curve, we find
\begin{align}\label{chiSchi2}
\chi_k \mathcal{S}_\Gamma \chi_i r^+ 
= \begin{cases}
\chi_k r^+, & k = i, \\
0, & k \neq i,
\end{cases}
\end{align}
and
\begin{align}\label{chiSchi}
\chi_k \mathcal{S}_\Gamma \chi_i r_j^- = \begin{cases}
-\chi_kr_j^-, & k = i = j, \\
\chi_k r_j^-, & k = i \neq j, \\
-2\chi_k r_j^-, & k \neq i = j, \\
0, & k \neq i \neq j,
\end{cases}
\end{align}
a.e. on $\Gamma$.
Equation (\ref{chiSchi}) implies 
\begin{align*}
\mathcal{S}_\Gamma r_j^-
& = \sum_{k,i = 1}^n \chi_k \mathcal{S}_\Gamma \chi_i r_j^-
= \chi_j \mathcal{S}_\Gamma \chi_j r_j^-
+ \sum_{k = i \neq j} \chi_k \mathcal{S}_\Gamma \chi_i r_j^-
+ \sum_{k \neq i = j} \chi_k\mathcal{S}_\Gamma \chi_i r_j^-
	\\
& = -\chi_jr_j^- + \sum_{k \neq j} \chi_k r_j^- - 2 \sum_{k \neq j} \chi_k r_j^- = -r_j^-, \qquad
1 \leq j \leq n. 
\end{align*}
Thus $\mathcal{S}_\Gamma r^- = - r^-$. Similarly, equation (\ref{chiSchi2}) implies $\mathcal{S}_\Gamma r^+ = r^+$.
This proves (\ref{SGammar}).

Equation (\ref{SGammar}) implies that $\mathcal{S}_\Gamma^2 r = r$ for every $r \in R(\Gamma)$. Since $R(\Gamma)$ is dense in $L^p(\Gamma)$ (see Lemma 9.14 in \cite{BK1997}), it follows that $\mathcal{S}_\Gamma^2h = h$ for every $h \in L^p(\Gamma)$.

To prove (\ref{jumpChD}), we note that part $(a)$ yields
\begin{align}\label{CrpCrm}
(\mathcal{C} r^+)(z) = \begin{cases} r^+(z), & z \in D_+, \\
0, & z \in D_-, 
\end{cases} \qquad
(\mathcal{C} r^-)(z) = \begin{cases} 0, & z \in D_+, \\
-r^-(z), & z \in D_-.
\end{cases}
\end{align}
It follows that $\mathcal{C}r \in \dot{E}^p(D_+ \cup D_-)$ for every $r \in R(\Gamma)$
and that 
\begin{align}\label{mathcalCr}
\mathcal{C}_+ r = r^+, \qquad \mathcal{C}_- r = -r^-.
\end{align}
Equations (\ref{SGammar}) and (\ref{mathcalCr}) imply
$$\frac{1}{2}(I + \mathcal{S}_\Gamma) r = 
\frac{1}{2}(r^+ + r^- + r^+ -r^-) = r^+ = \mathcal{C}_+ r .$$
Similarly,
$$\frac{1}{2}(-I + \mathcal{S}_\Gamma) r = \mathcal{C}_- r.$$
This shows that the Sokhotski-Plemelj formulas (\ref{plemelj}) are valid for all $r \in R(\Gamma)$. 

Let $h \in L^p(\Gamma)$. Let $r_n$ be a sequence in $R(\Gamma)$ converging to $h$ in $L^p(\Gamma)$. The boundedness of $\mathcal{S}_\Gamma$ on $L^p(\Gamma)$ implies
$$\mathcal{C}_\pm r_n = \frac{1}{2}(\pm r_n + \mathcal{S}_\Gamma r_n) \to \frac{1}{2}(\pm h + \mathcal{S}_\Gamma h) \quad \text{in}Ê\quad L^p(\Gamma)$$ 
as $n \to\infty$. Hence Lemma \ref{Epconvergencelemma} applied to each component of $\hat{\C} \setminus \Gamma$ implies that there exists a function $f \in \dot{E}^p(D_+ \cup D_-)$ such that $(\mathcal{C}r_n)|_{D_+ \cup D_-} \to f$ uniformly on compact subsets of $D_+ \cup D_-$ and $f_\pm = \frac{1}{2}(\pm h + \mathcal{S}_\Gamma h)$. Since $\mathcal{C}r_n \to \mathcal{C}h$ pointwise in $D_+$, we infer that $\mathcal{C}h = f \in \dot{E}^p(D_+ \cup D_-)$.
This proves (\ref{jumpChD}) in the case of $\infty \notin \Gamma$. It also follows that 
$$\mathcal{C}_\pm h = f_\pm = \frac{1}{2}(\pm h + \mathcal{S}_\Gamma h),$$
showing that the Sokhotski-Plemelj formulas (\ref{plemelj}) are valid for all $h \in L^p(\Gamma)$.

Suppose now that $\infty \in \Gamma$. Pick $z_0 \in D_-$ and let $\varphi(z) = 1/(z-z_0)$. 
Since $h \in \dot{L}^p(\Gamma)$, Lemma \ref{Lpdotlemma} implies that $\Phi h \in L^p(\varphi(\Gamma))$ and
\begin{align}\label{CGammahPsi}
(\mathcal{C}_\Gamma h)(z) = (\Psi^{-1} \mathcal{C}_{\varphi(\Gamma)} \Phi h)(z) \quad \text{for} \quad  z \in \C \setminus \Gamma.
\end{align}
The result of the previous paragraph implies that $\mathcal{C}_{\varphi(\Gamma)} \Phi h \in \dot{E}^p(\varphi(D_+ \cup D_-))$. Hence, in view of Proposition \ref{EpDinvariantprop} and equation (\ref{CGammahPsi}), $\mathcal{C}_\Gamma h \in \dot{E}^p(D_+ \cup D_-)$,
which proves (\ref{jumpChD}).
Similarly, the identity $\mathcal{S}_\Gamma^2 = I$ follows from the identity $\mathcal{S}_{\varphi(\Gamma)}^2 =I$ and equation (\ref{SGammainvariance}):
$$\mathcal{S}_\Gamma^2h = \Phi^{-1}\mathcal{S}_{\varphi(\Gamma)}^2 \Phi h = \Phi^{-1} \Phi h = h \quad \text{for all} \quad h \in \dot{L}^p(\Gamma).$$
This completes the proof of Theorem \ref{jumpth1}.

\subsection{Proof of Theorem \ref{jumpth2}}
We already established the Sokhotski-Plemelj formulas (\ref{plemelj}) in the case of $\infty \notin \Gamma$ (see the proof of part $(b)$ of Theorem \ref{jumpth1}). So suppose $\infty \in \Gamma$. Pick $z_0 \in D_-$ and let $\varphi(z) = 1/(z-z_0)$. The fact that $\infty \notin \varphi(\Gamma)$ together with the transformation properties (\ref{mathcalCinv}) and (\ref{SGammainvariance}) of $\mathcal{C}$ and $\mathcal{S}$ imply
$$\frac{1}{2}(\pm I + \mathcal{S}_\Gamma)h
= \frac{1}{2}\Phi^{-1}(\pm I + \mathcal{S}_{\varphi(\Gamma)})\Phi h
= (\Psi^{-1}\mathcal{C}_{\varphi(\Gamma)} \Phi h)_\pm
= (\mathcal{C}_\Gamma h)_\pm.$$
This completes the proof of (\ref{plemelj}).

The Sokhotski-Plemelj formulas (\ref{plemelj}) together with the fact that $\mathcal{S}_\Gamma^2 = I$ immediately imply that $\mathcal{C}_\pm$ are bounded orthogonal projections on $\dot{L}^p(\Gamma)$.

If $h = \mathcal{C}_+h - \mathcal{C}_-h \in \dot{L}^p(\Gamma)$, then $\mathcal{C} h \in \dot{E}^p(D_+ \cup D_-)$ by (\ref{jumpChD}). Hence equations (\ref{jumpCfplus}) and (\ref{jumpCfminus}) imply
$$(\mathcal{C} \mathcal{C}_+ h)(z)
= \begin{cases} (\mathcal{C}h)(z), & z \in D_+, \\
0, & z \in D_-,
\end{cases}
\qquad
(\mathcal{C} \mathcal{C}_- h)(z)
= \begin{cases} 0, & z \in D_+, \\
-(\mathcal{C}h)(z), & z \in D_-.
\end{cases}$$
These equations yield (\ref{ChDChD}). The last two statements of Theorem \ref{jumpth2} are easy consequences of (\ref{ChDChD}) and Theorem \ref{jumpth1}. This completes the proof. 
\proofend

\section{Riemann-Hilbert problems}\label{rhsec}\nequation
With Theorems \ref{jumpth1} and \ref{jumpth2} at our disposal, we can introduce a notion of $L^p$-RH problem for Carleson jump contours. Throughout this section, $\Gamma \subset \hat{\C}$ will denote a Carleson jump contour, $D_\pm \subset \hat{\C}$ will denote the associated open sets such that $\partial D_+ = -\partial D_- = \Gamma$, and we will assume that $1 < p <\infty$. We let $D = D_+ \cup D_-$.  

\subsection{Definition} \label{defsubsec}
Let $n \geq 1$ be an integer. Given an $n \times n$-matrix valued function $v: \Gamma \to GL(n, \C)$, we define a {\it solution of the $L^p$-RH problem determined by $(\Gamma, v)$} to be an $n \times n$-matrix valued function $m \in I + \dot{E}^p(D)$ such that the nontangential boundary values $m_\pm$ satisfy $m_+ = m_- v$ a.e. on $\Gamma$. 

\subsection{Properties of $m_\pm$}
In order to make contact with earlier works on $L^p$-RH problems on smooth contours, we show that $m$ is a solution of the $L^p$-RH problem if and only if the boundary functions $m_+$ and $m_-$ satisfy the properties (RH1)-(RH2) below. 

\begin{proposition}\label{equivalentdefprop}
Suppose $v: \Gamma \to GL(n, \C)$.
If $m \in I + \dot{E}^p(D)$ satisfies the $L^p$-RH problem determined by $(\Gamma, v)$, then the nontangential boundary values $m_\pm \in I + \dot{L}^p(\Gamma)$ satisfy the following two properties:
\begin{enumerate}[$(i)$]
\item[(RH1)] There exists a function $h \in \dot{L}^p(\Gamma)$ such that
\begin{align}\label{mpmCpmh}
m_\pm - I = \mathcal{C}_\pm h \quad \text{in}\quad \dot{L}^p(\Gamma).
\end{align}
\item[(RH2)] $m_+ = m_- v$ a.e. on $\Gamma$.
\end{enumerate}

Conversely, if $m_\pm \in I + \dot{L}^p(\Gamma)$ are a pair of $n\times n$-matrix valued functions satisfying (RH1) and (RH2), then $m = I + \mathcal{C}(m_+ - m_-) \in I + \dot{E}^p(D)$ satisfies $L^p$-RH problem determined by $(\Gamma, v)$.
\end{proposition}
\proofbegin
Theorem \ref{jumpth1} implies that if $m \in I + \dot{E}^p(D)$ satisfies the $L^p$-RH problem determined by $(\Gamma, v)$, then $m_\pm \in I + \dot{L}^p(\Gamma)$  and $m = I + \mathcal{C}(m_+ - m_-)$. Thus (RH1) is satisfied with $h = m_+ - m_-$. The property (RH2) holds by definition. 

Conversely, suppose $m_\pm \in I + \dot{L}^p(\Gamma)$ satisfy (RH1) and (RH2). By (RH1), $m_\pm \in I + \mathcal{C}_\pm  \dot{L}^p(\Gamma)$. Thus, Theorems \ref{jumpth1} and \ref{jumpth2} imply that $m_\pm$ are the nontangential boundary values of the function $m$ defined by $m = I + \mathcal{C}(m_+ - m_-) \in I + \dot{E}^p(D)$. It follows that $m$ satisfies the $L^p$-RH problem determined by $(\Gamma, v)$.
\proofend

\begin{remark}\upshape
In most earlier references on $L^p$-RH problems \cite{DZ2002a, DZ2002b, FIKN2006, Z1989}, a solution of an $L^p$-RH problem is defined as a pair of functions $m_\pm \in I + L^p(\Gamma)$ satisfying (RH1)-(RH2) (or properties very similar to (RH1)-(RH2)); the associated function $m(z)$ is then referred to as the `extension of $m_\pm$'. Here, in an effort to mimic the classical formulation of a RH problem as closely as possible, we have chosen to define a solution directly in terms of $m$. 
Proposition \ref{equivalentdefprop} shows that in the set-up provided by the spaces $\dot{L}^p(\Gamma)$ and $\dot{E}^p(D)$, the definitions in terms of $m$ and $m_\pm$ are equivalent.
\end{remark}

\begin{remark}\upshape
Condition (RH1) is equivalent to the condition that $m_\pm \in I + \mathcal{C}_\pm \dot{L}^p(\Gamma)$.
\end{remark}

\subsection{Uniqueness results}
We will show that the solution of the $L^p$-RH problem determined by $(\Gamma, v)$ is unique provided that $\det v = 1$ and $n \leq p$. 

\begin{lemma}\label{uniquelemma}
Suppose $v: \Gamma \to GL(n, \C)$.
Let $1 < p < \infty$ and define $q$ by $1/p + 1/q =1$. Let $m, \tilde{m} \in I + \dot{E}^p(D)$ be two solutions of the $L^p$-RH problem determined by $(\Gamma,v)$. If $m^{-1} \in I + \dot{E}^q(D)$, then $m(z) = \tilde{m}(z)$ for all $z \in D$.
\end{lemma}
\proofbegin  
Suppose $m, \tilde{m} \in I + \dot{E}^p(D)$ are two solutions of the $L^p$-RH problem determined by $(\Gamma,v)$. Suppose $m^{-1} \in I + \dot{E}^q(D)$. By Lemma \ref{Epdotlemma},
$$\tilde{m} m^{-1} - I = (\tilde{m} - I)(m^{-1} - I) + (\tilde{m} - I) 
+ (m^{-1} - I) \in \dot{E}^1(D) + \dot{E}^p(D) + \dot{E}^q(D) \subset \dot{E}^1(D).$$
Using Theorem \ref{jumpth1} and the fact that $(\tilde{m} m^{-1})_+ = \tilde{m}_{-}v v^{-1} m_{-}^{-1} = (\tilde{m} m^{-1})_-$ a.e on $\Gamma$, we find
$$\tilde{m} m^{-1} - I = \mathcal{C}((\tilde{m} m^{-1} - I)_+ - (\tilde{m} m^{-1} - I)_-)
= 0 \quad \text{on} \;\; D.$$
It follows that $m = \tilde{m}$ on $D$. 
\proofend

\begin{remark}\upshape
The assumption in Lemma \ref{uniquelemma} that $m^{-1} \in I + \dot{E}^q(D)$ implies that $m_\pm$ deliver a so-called $L^p$-canonical factorization of $v$; the uniqueness of the latter is known, see \cite{GKS2003, LS1987}.
\end{remark}

Suppose $v: \Gamma \to GL(2, \C)$ satisfies $\det v = 1$ a.e. on $\Gamma$.
If $m \in I + \dot{E}^2(D)$ is a solution of the $L^2$-RH problem determined by $(\Gamma,v)$, then Lemma \ref{Epdotlemma} shows that
$$\det m - 1 = (m_{11} - 1)(m_{22} - 1)  + (m_{11} - 1) + (m_{22} - 1) - m_{12} m_{21} \in \dot{E}^1(D).$$
By Theorem \ref{jumpth1} and the fact that $(\det m)_+ = (\det m)_-$ a.e. on $\Gamma$, we find
$$\det m - 1 = \mathcal{C}((\det m - 1)_+ - (\det m - 1)_-) = 0  \quad \text{on} \;\; D.$$
Thus, 
$$m^{-1} = \begin{pmatrix} m_{22} & -m_{12} \\ -m_{21} & m_{11} \end{pmatrix} \in I + \dot{E}^2(D).$$
Lemma \ref{uniquelemma} therefore implies that the solution of the $L^2$-RH problem determined by $(\Gamma,v)$ is unique if it exists.
This proves the special case $n = p = 2$ of the following theorem, which states that if $p \geq n$ and the $n\times n$-matrix valued jump function $v$ satisfies $\det v = 1$, then the solution of the $L^p$-RH problem determined by $(\Gamma,v)$ is unique if it exists.

Recall that the adjugate $\adj A$ of an $n \times n$ matrix $A$ is defined by 
$$(\adj A)_{ij} = (-1)^{i+j} m_{ji}(A), \qquad i,j = 1, \dots, n,$$
where $m_{ij}(A)$ denotes the $(ij)$th minor of $A$. By Cramer's rule, the inverse of $A$ is given by $A^{-1} = \adj(A)/\det(A)$ whenever $\det(A) \neq 0$. We continue to assume that $1 < p < \infty$.

\begin{theorem}\label{uniqueth}
Suppose $1 \leq n \leq p$. Suppose $v: \Gamma \to GL(n, \C)$ satisfies $\det v = 1$ a.e. on $\Gamma$.
\begin{enumerate}[$(a)$]
\item If $m$ is a solution of the $L^p$-RH problem determined by $(\Gamma,v)$, then $\det m(z) = 1$ for all $z \in D$.

\item The solution of the $L^p$-RH problem determined by $(\Gamma,v)$ is unique if it exists.
\end{enumerate}
\end{theorem}
\proofbegin 
Let $m \in I + \dot{E}^p(D)$ be a solution of the $L^p$-RH problem determined by $(\Gamma,v)$ for some $p \geq n$.
By Lemma \ref{Epdotlemma}, if $\{f_j\}_1^k \subset \dot{E}^p(D)$ and $1 \leq k \leq n$, then $\Pi_{j=1}^k f_j \in \dot{E}^{p/k}(D) \subset \dot{E}^1(D)$. 
It follows that
$$\det m - 1 = \det(I + (m- I)) - 1 \in \dot{E}^1(D).$$
Using Theorem \ref{jumpth1} and the fact that $(\det m)_+ = (\det m)_-$ a.e. on $\Gamma$, we find
$$\det m - 1 = \mathcal{C}((\det m)_+ - (\det m)_-) = 0  \quad \text{on} \;\; D,$$
which proves $(a)$. To prove $(b)$, we note that Lemma \ref{Epdotlemma} implies 
$$\adj m \in I + \dot{E}^{\frac{p}{n-1}}(D) \subset  I + \dot{E}^{\frac{p}{p-1}}(D).$$ 
Hence, by $(a)$,
$$m^{-1} = \adj m \in I + \dot{E}^{\frac{p}{p-1}}(D),$$
so that $(b)$ follows from Lemma \ref{uniquelemma}.
\proofend

\begin{remark}\upshape
For a sufficiently smooth contour, the special case $n = p = 2$ of Theorem \ref{uniqueth} was proved in \cite{D1999, DZ2003}. Theorem \ref{uniqueth} generalizes this result to the case of a Carleson contour $\Gamma$ and any $1 \leq n \leq p$. As an application, we note that the case $n=3$ is relevant for the $3 \times 3$-matrix RH problem associated with the Degasperis-Procesi equation, see Figure \ref{Dns.pdf}.
\end{remark}

\subsection{A singular integral equation}
Given a Banach space $X$, we let $\mathcal{B}(X)$ denote the space of bounded linear operators on $X$. 
Given two function $w^\pm \in \dot{L}^p(\Gamma) \cap L^\infty(\Gamma)$, we define the operator $\mathcal{C}_{w}: \dot{L}^p(\Gamma) + L^\infty(\Gamma) \to \dot{L}^p(\Gamma)$ by 
$$\mathcal{C}_{w}(f) = \mathcal{C}_+(f w^-) + \mathcal{C}_-(f w^+).$$
We fix a point $z_0 \in \C \setminus \Gamma$ and let $\| \cdot \|_{\dot{L}^p(\Gamma)}$ denote the associated norm on $\dot{L}^p(\Gamma)$ defined in (\ref{dotLpnorm}).
The estimate
\begin{align*}
\|\mathcal{C}_wf\|_{\dot{L}^p(\Gamma)} 
& = \|\mathcal{C}_+(f w^-) + \mathcal{C}_-(f w^+)\|_{\dot{L}^p(\Gamma)} 
	\\
& \leq C \|f\|_{\dot{L}^p(\Gamma)} \max\big\{\|w^+\|_{L^\infty(\Gamma)}, \|w^-\|_{L^\infty(\Gamma)} \big\} \quad \text{for} \quad f \in \dot{L}^p(\Gamma),
\end{align*}
where $C = \max\{\|\mathcal{C}_+\|_{\mathcal{B}(\dot{L}^p(\Gamma))}, \|\mathcal{C}_-\|_{\mathcal{B}(\dot{L}^p(\Gamma))}\} < \infty$, implies that
\begin{align}\label{Cwnorm}
\|\mathcal{C}_w\|_{\mathcal{B}(\dot{L}^p(\Gamma))} \leq C \max\big\{\|w^+\|_{L^\infty(\Gamma)}, \|w^-\|_{L^\infty(\Gamma)} \big\}.
\end{align}

The next proposition shows that if $v = (v^-)^{-1}v^+$ and $w^\pm = \pm v^\pm \mp I$ then the $L^p$-RH problem determined by $(\Gamma, v)$ is equivalent to the following singular integral equation for $\mu \in I + \dot{L}^p(\Gamma)$ cf. \cite{BC1984}:
\begin{align}\label{rhoeq}
\mu - I = \mathcal{C}_w(\mu)  \quad \text{in}\quad \dot{L}^p(\Gamma).
\end{align}

\begin{proposition}\label{muprop}
Given $v^\pm: \Gamma \to GL(n, \C)$, let $v = (v^-)^{-1}v^+$, $w^+ = v^+ - I$, and $w^- = I - v^-$. Suppose $v^\pm, (v^\pm)^{-1} \in I +  \dot{L}^p(\Gamma) \cap L^\infty(\Gamma)$.
If $m \in I + \dot{E}^p(D)$ satisfies the $L^p$-RH problem determined by $(\Gamma, v)$, then $\mu = m_+ (v^+)^{-1} = m_- (v^-)^{-1} \in I + \dot{L}^p(\Gamma)$ satisfies (\ref{rhoeq}). 
Conversely, if $\mu \in I + \dot{L}^p(\Gamma)$ satisfies (\ref{rhoeq}), then
$m = I + \mathcal{C}(\mu(w^+ + w^-)) \in I + \dot{E}^p(D)$ satisfies the $L^p$-RH problem determined by $(\Gamma, v)$. 
\end{proposition}
\proofbegin
Suppose $m \in I + \dot{E}^p(D)$ satisfies the $L^p$-RH problem determined by $(\Gamma, v)$ and let $\mu = m_+ (v^+)^{-1} = m_- (v^-)^{-1}$.
By Theorem \ref{jumpth1}, $m_\pm = I + \dot{L}^p(\Gamma)$ and hence $\mu \in I + \dot{L}^p(\Gamma)$. Moreover, by Theorem \ref{jumpth2},
\begin{align*}
\mathcal{C}_w\mu & = \mathcal{C}_+(\mu(I - v^-) ) - \mathcal{C}_-(\mu (I - v^+))
	\\
&= \mathcal{C}_+(\mu - I + I - m_-) - \mathcal{C}_-(\mu - I + I - m^+))
 = (\mathcal{C}_+ - \mathcal{C}_-)(\mu - I) = \mu - I.
\end{align*}

Conversely, suppose $\mu \in I + \dot{L}^p(\Gamma)$ satisfies (\ref{rhoeq}). The assumption $v^\pm \in I + \dot{L}^p(\Gamma) \cap L^\infty(\Gamma)$ implies that $\mu w^\pm \in \dot{L}^p(\Gamma)$. Hence $m = I + \mathcal{C} (\mu (w^+ + w^-)) \in I + \dot{E}^p(D)$ and
\begin{subequations}\label{mCCmuw}
\begin{align}
& m_+ = I + (\mathcal{C}_+ - \mathcal{C}_-)(\mu w^+) + \mathcal{C}_w\mu
  = \mu (w^+ + I) = \mu v^+,
  	\\
& m_- = I - (\mathcal{C}_+ - \mathcal{C}_-)(\mu w^-) + \mathcal{C}_w\mu
  = \mu (I - w^-) = \mu v^-.
\end{align}
\end{subequations}
It follows that $m_+ = m_- v$ a.e. on $\Gamma$.
\proofend

\subsection{Vanishing lemma}
Given $v: \Gamma \to GL(n, \C)$, we define a {\it solution of the homogeneous $L^p$-RH problem determined by $(\Gamma, v)$} to be an $n \times n$-matrix valued function $m \in \dot{E}^p(D)$ such that $m_+ = m_- v$ a.e. on $\Gamma$. 

\begin{lemma}\label{abcdlemma}
Given $v^\pm: \Gamma \to GL(n, \C)$, let $v = (v^-)^{-1}v^+$, $w^+ = v^+ - I$, and $w^- = I - v^-$. Suppose $v^\pm, (v^\pm)^{-1} \in I +  \dot{L}^p(\Gamma) \cap L^\infty(\Gamma)$.
Then the implications
$$(a) \implies (b) \implies (c) \implies (d)$$
are valid for the following statements:
\begin{enumerate}[$(a)$]
\item The map $I - \mathcal{C}_w: \dot{L}^p(\Gamma)Ê\to \dot{L}^p(\Gamma)$ is bijective.

\item The $L^p$-RH problem determined by $(\Gamma, v)$ has a unique solution. 

\item The homogeneous $L^p$-RH problem determined by $(\Gamma, v)$ has only the zero solution.

\item  The map $I - \mathcal{C}_w: \dot{L}^p(\Gamma)Ê\to \dot{L}^p(\Gamma)$ is injective.
\end{enumerate}
\end{lemma}
\proofbegin
$(a) \implies (b)$ \; Suppose $I - \mathcal{C}_w:\dot{L}^p(\Gamma)Ê\to \dot{L}^p(\Gamma)$ is a bijection. Then 
$\mu = I + (I - \mathcal{C}_w)^{-1}\mathcal{C}_wI \in I + \dot{L}^p(\Gamma)$
is a solution of (\ref{rhoeq}). Hence, by Proposition \ref{muprop}, $m = I + \mathcal{C}(\mu (w^+ + w^-)) \in I + \dot{E}^p(D)$ satisfies the $L^p$-RH problem determined by $(\Gamma, v)$. Moreover, by (\ref{mCCmuw}), $m_\pm = \mu v^\pm$.
If $\tilde{m}Ê\in I + \dot{E}^p(D)$ is another solution of this RH problem, then Proposition \ref{muprop} implies that $\tilde{\mu} = \tilde{m}_\pm (v^\pm)^{-1} \in I + \dot{L}^p(\Gamma)$ is a solution of (\ref{rhoeq}). But then $ \tilde{\mu} = I + (I - \mathcal{C}_w)^{-1}\mathcal{C}_wI = \mu$, so that $\tilde{m}_\pm = m_\pm$ a.e. on $\Gamma$. Theorem \ref{jumpth1} now yields
$$m = I + \mathcal{C}(m_+ - m_-) = I + \mathcal{C}(\tilde{m}_+ - \tilde{m}_-) = \tilde{m} \quad \text{in} \;\; \dot{E}^p(D),$$
showing that the solution is unique.

$(b) \implies (c)$ \; Let $mÊ\in I + \dot{E}^p(D)$ be the unique solution of the $L^p$-RH problem determined by $(\Gamma, v)$ and suppose $\tilde{m} \in \dot{E}^p(D)$ satisfies the homogeneous RH problem determined by $(\Gamma, v)$. By Proposition \ref{muprop}, $\mu = m_- (v^-)^{-1}$ satisfies equation (\ref{rhoeq}). By Theorem \ref{jumpth2}, 
\begin{align*}
(I - \mathcal{C}_w)(\tilde{m}_-(v^-)^{-1}) 
& = \tilde{m}_-(v^-)^{-1} - \mathcal{C}_+(\tilde{m}_- (v^-)^{-1}(I - v^-)) - \mathcal{C}_-(\tilde{m}_- (v^-)^{-1} (v^+ - I))
	\\
& = \tilde{m}_-(v^-)^{-1} - \mathcal{C}_+(\tilde{m}_- (v^-)^{-1} - \tilde{m}_-) 
- \mathcal{C}_-(\tilde{m}_+ - \tilde{m}_- (v^-)^{-1})
= 0.
\end{align*}
Hence $\mu = (m_- + \tilde{m}_-)(v^-)^{-1}$ also satisfies equation (\ref{rhoeq}).
By Proposition \ref{muprop} and uniqueness of $m$, we conclude that
$$m = I + \mathcal{C}(m_-(v^-)^{-1} (w^+ + w^-)) = I + \mathcal{C}((m_- + \tilde{m}_-) (v^-)^{-1} (w^+ + w^-)).$$
But then $\mathcal{C}_\pm(\tilde{m}_- (v^-)^{-1} (w^+ + w^-)) = 0$, and so 
$$\tilde{m}_+ - \tilde{m}_- = \tilde{m}_- (v^-)^{-1} (w^+ + w^-) = (\mathcal{C}_+ - \mathcal{C}_-)(\tilde{m}_- (v^-)^{-1}  (w^+ + w^-)) = 0$$ 
a.e. on $\Gamma$. Thus, by Theorem \ref{jumpth1}, $\tilde{m} = \mathcal{C}(\tilde{m}_+ - \tilde{m}_-) = 0$.

$(c) \implies (d)$ \; 
Suppose the homogeneous $L^p$-RH problem determined by $(\Gamma, v)$ has only the zero solution. Suppose $h \in \dot{L}^p(\Gamma)$ satisfies $(I - \mathcal{C}_w)h = 0$. Let $m = \mathcal{C}(h (w^+ + w^-)) \in \dot{E}^p(D)$. Since 
$$m_+ = \mathcal{C}_+ (h (w^+ + w^-)) 
= \mathcal{C}_+ (h w^+) + h -  \mathcal{C}_- (h w^+)
= h w^+ + h = h v^+,$$
$$m_- = \mathcal{C}_- (h (w^+ + w^-)) 
= h -  \mathcal{C}_+ (h w^-) + \mathcal{C}_- (h w^-)
= h - hw^-
= hv^-,$$
it follows that $m_+ = m_-v$ a.e. on $\Gamma$. Hence $m = 0$ by uniqueness of the solution of the homogeneous problem. Thus $h = m_-(v^-)^{-1} = 0$, showing that $(I - \mathcal{C}_w)$ is injective.
\proofend

Let $C(\Gamma)$ denote the set of restrictions to $\Gamma$ of continuous functions $\hat{\C} \to \C$. If $\Gamma \subset \hat{\C}$ is given the subspace topology, Tietze's extension theorem implies that $C(\Gamma)$ coincides with the set of continuous functions $\Gamma \to \C$.
We will show that if $w^\pm \in C(\Gamma)$ then the operator $I - \mathcal{C}_w$ is Fredholm. If, in addition, $w^\pm$ are nilpotent, the Fredholm index of this operator is zero, so that all four statements $(a)$-$(d)$ of Lemma \ref{abcdlemma} are equivalent. 

For a Banach space $X$, let $\mathcal{K}(X) \subset \mathcal{B}(X)$ denote the set of compact operators on $X$. The set of Fredholm operators $\mathcal{F}(X)$ is open in $\mathcal{B}(X)$ and the index map $\Ind:\mathcal{F}(X) \to \Z$ is constant on the connected components of $\mathcal{F}(X)$. 
If $X = \dot{L}^p(\Gamma)$, we define $\mathcal{B}(X)$, $\mathcal{K}(X)$, and $\mathcal{F}(X)$ as the set of bounded, compact, and Fredholm operators on $L^p(\Gamma, w)$ where $w(z) =  |z - z_0|^{1 - \frac{2}{p}}$ and $z_0$ is any point of $\C \setminus \Gamma$. 

Given $w^\pm, \tilde{w}^\pm \in L^\infty(\Gamma)$ such that $\tilde{w}^+ = (w^+ + I)^{-1} - I$ and $\tilde{w}^- = I - (I - w^-)^{-1}$, we define $T_w, T_{\tilde{w}}:\dot{L}^p(\Gamma)Ê\to \dot{L}^p(\Gamma)$ by
\begin{subequations}\label{Twdef} 
\begin{align}
& T_w = \mathcal{C}_+ \mathcal{R}_{\tilde{w}^-} \mathcal{C}_- \mathcal{R}_{w^+ + w^-}
+ \mathcal{C}_- \mathcal{R}_{\tilde{w}^+} \mathcal{C}_+ \mathcal{R}_{w^+ + w^-},
	\\
& T_{\tilde{w}} = \mathcal{C}_+ \mathcal{R}_{w^-} \mathcal{C}_- \mathcal{R}_{\tilde{w}^+ + \tilde{w}^-}
+ \mathcal{C}_- \mathcal{R}_{w^+} \mathcal{C}_+ \mathcal{R}_{\tilde{w}^+ + \tilde{w}^-},
\end{align}
\end{subequations}
where the right multiplication operator $\mathcal{R}_g$ is defined for functions $g(z)$ and $h(z)$ by 
$$(\mathcal{R}_g h)(z)= h(z) g(z).$$


\begin{theorem}
Given $v^\pm: \Gamma \to GL(n, \C)$, let $v = (v^-)^{-1}v^+$, $w^+ = v^+ - I$, and $w^- = I - v^-$. Suppose $v^\pm, (v^\pm)^{-1} \in  I + \dot{L}^p(\Gamma) \cap L^\infty(\Gamma)$ and $v^\pm \in C(\Gamma)$.
\begin{enumerate}[$(a)$]
\item The operator $I - \mathcal{C}_w:\dot{L}^p(\Gamma)Ê\to \dot{L}^p(\Gamma)$ is Fredholm.
\item If $w^\pm$ are nilpotent matrices, then $I - \mathcal{C}_w$ has Fredholm index zero; in this case, each of the four statements $(a)$-$(d)$ of Lemma \ref{abcdlemma} implies the other three.
\end{enumerate}
\end{theorem}
\proofbegin
Since $\Gamma \subset \hat{\C}$ is compact, there exists a $c$ such that $|\det v^\pm| \geq c > 0$ on $\Gamma$. Thus $(v^\pm)^{-1} \in C(\Gamma)$.
Let $\tilde{w}^+ = (v^+)^{-1} - I$ and $\tilde{w}^- = I - (v^-)^{-1}$. Then $\mathcal{C}_{w}$ and $\mathcal{C}_{\tilde{w}}$ are bounded $\dot{L}^p(\Gamma)Ê\to \dot{L}^p(\Gamma)$. 

Assume first that $\infty \notin \Gamma$. 

{\it Step 1.} We will show that $T_w$ and $T_{\tilde{w}}$ defined by (\ref{Twdef}) are compact operators on $L^p(\Gamma)$. By Mergelyan's rational approximation theorem (see p. 119 of \cite{G1987}), $R(\Gamma)$ is dense in $C(\Gamma)$ equipped with the $L^\infty$-norm. Let $\{w_n^\pm\}_1^\infty \subset R(\Gamma)$ be sequences such that $\lim_{n\to \infty} \|w^\pm - w_n^\pm\|_{L^\infty(\Gamma)} = 0$.
Since
$$((\mathcal{R}_{w_n^+} \mathcal{S}_\Gamma - \mathcal{S}_\Gamma \mathcal{R}_{w_n^+})h)(z) 
= \frac{1}{\pi i}Ê\int_{\Gamma} h(z') \frac{w_n^+(z) - w_n^+(z')}{z' - z} dz', \qquad z \in \Gamma,$$
the operators $\mathcal{R}_{w_n^+} \mathcal{S}_\Gamma - \mathcal{S}_\Gamma \mathcal{R}_{w_n^+}$ are integral operators with continuous kernels. A standard argument based on Ascoli's theorem implies that they are compact $L^p(\Gamma) \to C(\Gamma)$; hence they are also compact $L^p(\Gamma) \to L^p(\Gamma)$. Since
$$\|(\mathcal{R}_{w^+} \mathcal{S}_\Gamma - \mathcal{S}_\Gamma \mathcal{R}_{w^+}) - (\mathcal{R}_{w_n^+} \mathcal{S}_\Gamma - \mathcal{S}_\Gamma \mathcal{R}_{w_n^+})\|_{\mathcal{B}(L^p(\Gamma))}
\leq 2 \| w^+ - w_n^+\|_{L^\infty(\Gamma)} \|\mathcal{S}_\Gamma \|_{\mathcal{B}(L^p(\Gamma))} \to 0$$
as $n \to \infty$, it follows that 
$$\mathcal{R}_{w^+} \mathcal{S}_\Gamma - \mathcal{S}_\Gamma \mathcal{R}_{w^+}
= 2(\mathcal{R}_{w^+} \mathcal{C}_+ - \mathcal{C}_+ \mathcal{R}_{w^+})$$
 is compact.
Since the compact operators form a two-sided ideal, we find that
$$\mathcal{C}_-(\mathcal{R}_{w^+} \mathcal{C}_+ - \mathcal{C}_+ \mathcal{R}_{w}^+)\mathcal{R}_{\tilde{w}^+ + \tilde{w}^-}
= \mathcal{C}_- \mathcal{R}_{w^+} \mathcal{C}_+ \mathcal{R}_{\tilde{w}^+ + \tilde{w}^-}$$
is a compact operator on $L^p(\Gamma)$. Similar arguments apply to the other terms in (\ref{Twdef}). This shows that $T_{w}$ and $T_{\tilde{w}}$ are compact on $L^p(\Gamma)$.

{\it Step 2.} We will show that $I - \mathcal{C}_w$ is Fredholm on $\dot{L}^p(\Gamma)$. 
Let $h \in L^p(\Gamma)$. Then
\begin{align*}
  \mathcal{C}_{\tilde{w}}\mathcal{C}_w &h
  = \mathcal{C}_+((\mathcal{C}_+(hw^-) + \mathcal{C}_-(hw^+))\tilde{w}^-)
 +  \mathcal{C}_-((\mathcal{C}_+(hw^-) + \mathcal{C}_-(hw^+))\tilde{w}^+)
	\\
 = &\; \mathcal{C}_+((hw^- + \mathcal{C}_-(hw^-) + \mathcal{C}_-(hw^+))\tilde{w}^-)
 +  \mathcal{C}_-((\mathcal{C}_+(hw^-) - hw^+ + \mathcal{C}_+(hw^+))\tilde{w}^+)
	\\
= &\; T_wh + \mathcal{C}_+(hw^-\tilde{w}^-) - \mathcal{C}_-(hw^+\tilde{w}^+).
\end{align*}
In view of the identities $w^+ \tilde{w}^+Ê= - w^+ - \tilde{w}^+$ and $w^- \tilde{w}^-Ê= w^- + \tilde{w}^-$, the right-hand side equals $T_wh + \mathcal{C}_wh + \mathcal{C}_{\tilde{w}} h$.
Hence
$$I + T_w = (I - \mathcal{C}_{\tilde{w}})(I - \mathcal{C}_w).$$
Interchanging $w$ and $\tilde{w}$ in the above argument, we find
$$I + T_{\tilde{w}} = (I - \mathcal{C}_w)(I - \mathcal{C}_{\tilde{w}}).$$
It follows that $I - \mathcal{C}_w$ is invertible modulo compact operators; hence $I - \mathcal{C}_w$ is Fredholm on $L^p(\Gamma)$.
Since the norms of $L^p(\Gamma)$ and $\dot{L}^p(\Gamma)$ are equivalent when $\Gamma$ is bounded, this proves $(a)$ in the case of $\infty \notin \Gamma$. 

{\it Step 3.} The map $t \mapsto I - \mathcal{C}_{tw}$ is continuousÊ $[0,1] \to \mathcal{B}(L^p(\Gamma))$ because
$$\| \mathcal{C}_{tw} -  \mathcal{C}_{sw} \|_{\mathcal{B}(L^p(\Gamma))}
= |t - s| \|\mathcal{C}_{w} \|_{\mathcal{B}(L^p(\Gamma))}, \qquad t,s \in [0,1].$$
If $w^\pm$ are nilpotent, then $tw^\pm \in C(\Gamma)$ and $\det(tw^+ + I) = \det(I - tw^-) = 1$, thus the operator $I - \mathcal{C}_{tw}$ is Fredholm on $L^p(\Gamma)$ for $t \in [0,1]$ by Step 2. 
Since the Fredholm index is constant on connected components, 
this proves $(b)$ in the case of $\infty \notin \Gamma$. 

{\it Step 4.} Suppose now that $\infty \in \Gamma$. Pick $z_0 \in D_-$. Let $\Phi:\dot{L}^p(\Gamma) \to L^p(\varphi(\Gamma))$ be the bijection defined in (\ref{Phidef}). Equipping $\dot{L}^p(\Gamma)$ with the norm (\ref{dotLpnorm}), $\Phi$ is an isometry by Lemma \ref{Lpdotlemma}. 
Let  $\mathcal{C} = \mathcal{C}_\Gamma$Ê and $\tilde{\mathcal{C}} = \mathcal{C}_{\varphi(\Gamma)}$ denote the Cauchy operators associated with the contours Ê$\Gamma$ and $\varphi(\Gamma)$, respectively. Using (\ref{mathcalCinv}), we find
\begin{align} \nonumber
I - \mathcal{C}_w & = I  - \mathcal{C}_+R_{w^-} - \mathcal{C}_-R_{w^+}
= \Phi^{-1}(I - \Phi \mathcal{C}_+ \Phi^{-1} \Phi R_{w^-} \Phi^{-1}
 - \Phi \mathcal{C}_- \Phi^{-1} \Phi R_{w^+} \Phi^{-1}) \Phi
	\\ \label{CwCGammaRw}
& = \Phi^{-1}(I -  \tilde{\mathcal{C}}_+ R_{w^- \circ \varphi^{-1}} - \tilde{\mathcal{C}}_- R_{w^+ \circ \varphi^{-1}}) \Phi
= \Phi^{-1}(I -  \tilde{\mathcal{C}}_{w \circ \varphi^{-1}}) \Phi.
\end{align}
Since $v^\pm \circ \varphi^{-1}: \varphi(\Gamma) \to GL(n, \C)$ satisfy $v^\pm \circ \varphi^{-1} \in C(\varphi(\Gamma))$ as well as 
$$v^\pm \circ \varphi^{-1}, (v^\pm)^{-1} \circ \varphi^{-1} \in I + L^p(\varphi(\Gamma)) \cap L^\infty(\varphi(\Gamma)),$$ 
Step 2 implies that the operator 
$I -  \tilde{\mathcal{C}}_{w \circ \varphi^{-1}} \in \mathcal{B}(L^p(\varphi(\Gamma)))$
is Fredholm. Since $\Phi$ is an isometry, equation (\ref{CwCGammaRw}) implies that $I - \mathcal{C}_w \in \mathcal{B}(\dot{L}^p(\Gamma))$ is also Fredholm of the same index.
\proofend

\subsection{Reversal of subcontours}\label{reversalsubsec}
It is sometimes convenient to consider RH problems with jumps across contours which are not Carleson jump contours but which can be turned into Carleson jump contours by reorienting an appropriate subcontour. We make the following definition: If $\tilde{\Gamma}$ denotes the Carleson jump contour $\Gamma$ with the orientation reversed on a subset $\Gamma_0 \subset \Gamma$ and $\tilde{v}$ is defined by
$$\tilde{v} = \begin{cases} v & \text{on $\Gamma \setminus \Gamma_0$}, \\
v^{-1} & \text{on $\Gamma_0$},
\end{cases}$$
then we say that $m \in I + \dot{E}^p(D)$ satisfies the $L^p$-RH problem determined by $(\tilde{\Gamma}, \tilde{v})$ if and only if $m$ satisfies the $L^p$-RH problem determined by $(\Gamma, v)$.

\subsection{Contour deformations}
Many applications of RH problems rely on arguments involving contour deformations. For example, in the nonlinear steepest descent method of \cite{DZ1993}, the jump contour is deformed in such a way that $w = v - I$ is exponentially small away from a finite number of critical points. Theorem \ref{deformationth} below gives conditions under which the deformed RH problem is equivalent to the original one.

\begin{lemma}\label{fglemma}
Let $D$ be the union of any number of components of $\hat{\C} \setminus \Gamma$, where $\Gamma$ is a Carleson jump contour. Let $E^\infty(D)$ denote the space of bounded analytic functions in $D$. 
If $f \in \dot{E}^p(D)$ and $g \in E^\infty(D)$, then $fg \in \dot{E}^p(D)$.
\end{lemma}
\proofbegin
The result is immediate when $\infty \notin \Gamma$. The case of $\infty \in \Gamma$ can be reduced to the case of $\infty \notin \Gamma$ by means of Proposition \ref{EpDinvariantprop}.
\proofend

\begin{figure}
\begin{center}
\quad  \begin{overpic}[width=.43\textwidth]{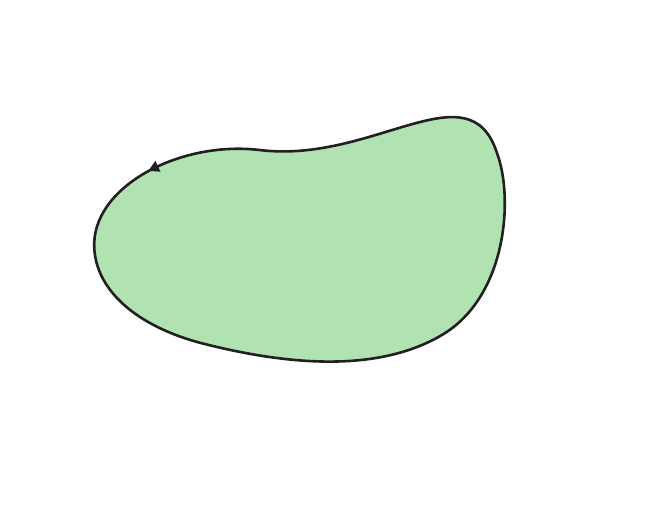}
 \put(50,40){$B_+$}
 \put(50,68){$B_-$}
 \put(22,56){$\gamma$}
 \end{overpic}
 \qquad
 \begin{overpic}[width=.45\textwidth]{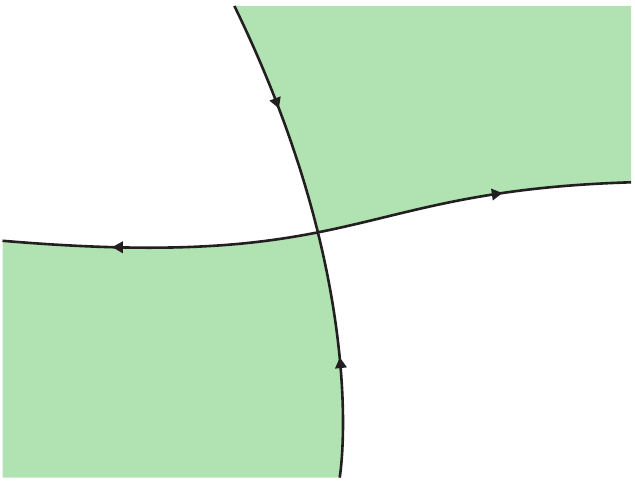}
 \put(70,60){$D_+$}
 \put(20,55){$D_-$}
 \put(20,16){$D_+$}
 \put(70,20){$D_-$}
  \put(77,38){$\Gamma$}
 \end{overpic}
    \bigskip
   \begin{figuretext}\label{gammafig1}
      The contours $\gamma$ and $\Gamma$.
      \end{figuretext}
   \end{center}
\end{figure}
 
 \begin{figure}
\begin{center}
\quad
 \begin{overpic}[width=.45\textwidth]{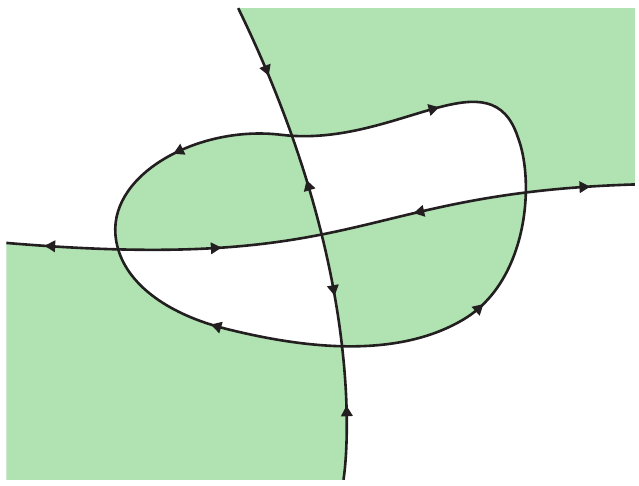}
 \put(85,62){$\hat{D}_+$}
 \put(15,60){$\hat{D}_-$}
 \put(15,11){$\hat{D}_+$}
 \put(80,11){$\hat{D}_-$}
 \put(30,43){$\hat{D}_+$}
 \put(36,29){$\hat{D}_-$}
 \put(62,31){$\hat{D}_+$}
 \put(60,48){$\hat{D}_-$}
   \put(6,40){$\hat{\Gamma}$}
 \end{overpic}
 \qquad
 \begin{overpic}[width=.45\textwidth]{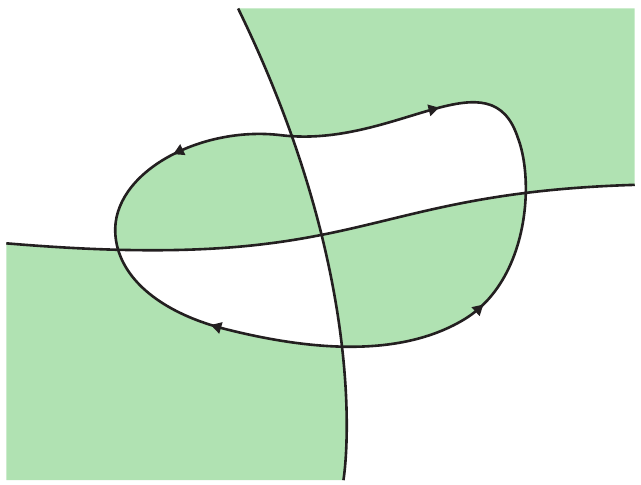}
 \put(30,43){$U_+$}
 \put(36,29){$U_-$}
 \put(62,31){$U_+$}
 \put(60,48){$U_-$}
 \put(23,56){$\gamma_+$}
 \put(30,19){$\gamma_-$}
 \put(75,22){$\gamma_+$}
 \put(65,63){$\gamma_-$}
  \end{overpic}
   \bigskip
   \begin{figuretext}\label{gammafig2}
      The contours $\hat{\Gamma} = \Gamma \cup \gamma$ and $\gamma_\pm$.     
    \end{figuretext}
   \end{center}
\end{figure}
Let $\hat{\Gamma} = \Gamma \cup \gamma$ denote the union of the Carleson jump contour $\Gamma$ and a curve $\gamma \in \mathcal{J}$, see Figures \ref{gammafig1} and \ref{gammafig2}. Suppose that, reversing the orientation on a subcontour if necessary, $\hat{\Gamma}$ is a Carleson jump contour. To be definite, we henceforth fix an orientation on the contour $\hat{\Gamma}$ which turns it into a Carleson jump contour, and we endow the contours $\Gamma$ and $\gamma$ with the orientations they inherit as subsets of $\hat{\Gamma}$. Then $\Gamma$ is a Carleson jump contour up to reorientation; we define a solution of the $L^p$-RH problem determined by $(\Gamma, v)$ as in Subsection \ref{reversalsubsec}. 

Let $B_+$ and $B_-$ denote the two components of $\hat{\C} \setminus \gamma$. Without loss of generality, we may assume that $\infty \in \bar{B}_-$.
Let $\hat{D}_\pm$ be the open sets such that $\hat{\C} \setminus \hat{\Gamma} = \hat{D}_+ \cup \hat{D}_-$ and $\partial \hat{D}_+ = - \partial \hat{D}_- = \hat{\Gamma}$. 
Let $U_\pm = \hat{D}_\pm \cap B_+$. Let $\hat{D} = \hat{D}_+ \cup \hat{D}_-$ and $U = U_+ \cup U_-$.
Let $\gamma_+$ and $\gamma_-$ be the parts of $\gamma$ that belong to the boundary of $U_+$ and $U_-$, respectively.  The orientations of $\gamma_\pm$ are such that $B_+$ lies to the left of $\gamma_+$, whereas $B_+$ lies to the right of $\gamma_-$.

\begin{theorem}\label{deformationth}
Let $1 < p < \infty$ and let $n \geq 1$ be an integer. Suppose $v: \Gamma \to GL(n, \C)$.
Suppose $m_0:U \to GL(n,\C)$ satisfies
$$m_0, m_0^{-1} \in I + \dot{E}^p(U) \cap E^\infty(U).$$
Define $\hat{v}:\hat{\Gamma} \to GL(n, \C)$ by
\begin{align*}
\hat{v} 
=  \begin{cases}
 m_{0-} v m_{0+}^{-1} & \text{on}Ê\quad  \Gamma \cap B_+, \\
m_{0+}^{-1} & \text{on}Ê\quad \gamma_+, \\
m_{0-} & \text{on}Ê\quad \gamma_-, \\
v & \text{on}Ê\quad \Gamma \cap B_-.
\end{cases}
\end{align*}

Then the $L^p$-RH problems determined by $(\Gamma,v)$ and $(\hat{\Gamma}, \hat{v})$ are equivalent in the following sense:
If $m \in I + \dot{E}^p(D)$ satisfies the $L^p$-RH problem determined by $(\Gamma,v)$, then the function $\hat{m}(z)$ defined for $z \in \hat{D}$ by
\begin{align}\label{hatmdefmm0}
\hat{m} = \begin{cases}
mm_0^{-1} & \text{on} \quad \hat{D} \cap B_+,\\
m & \text{on} \quad \hat{D} \cap B_-,
\end{cases}
\end{align}
satisfies the $L^p$-RH problem determined by $(\hat{\Gamma}, \hat{v})$.

Conversely, if $\hat{m} \in I + \dot{E}^p(\hat{D})$ satisfies the $L^p$-RH problem determined by $(\hat{\Gamma}, \hat{v})$, then the function $m(z)$ defined for $z \in \hat{D}$ by
\begin{align}\label{Amdef}
m = \begin{cases}
\hat{m}m_0 & \text{on} \quad \hat{D} \cap B_+,\\
\hat{m} & \text{on} \quad  \hat{D} \cap B_-,
\end{cases}
\end{align}
and extended to $D \cap \gamma$ by continuity, satisfies the $L^p$-RH problem determined by $(\Gamma, v)$.
\end{theorem}
\proofbegin
Suppose $m \in I + \dot{E}^p(D)$ satisfies the $L^p$-RH problem determined by $(\Gamma,v)$. Define $\hat{m}(z)$ for $z \in \hat{D}$ by (\ref{hatmdefmm0}). Using the identity $mm_0^{-1} = (m-I)m_0^{-1} + m_0^{-1} $ and Lemma \ref{fglemma}, we find that $\hat{m} \in I + \dot{E}^p(\hat{D})$.
The nontangential boundary values $\hat{m}_\pm \in I + \dot{L}^p(\hat{\Gamma})$ satisfy
$$\hat{m}_\pm = \begin{cases}
m_\pm m_{0\pm}^{-1} & \text{on} \quad \hat{\Gamma} \cap B_+,\\
m_\pm, & \text{on} \quad \hat{\Gamma} \cap B_-.
\end{cases}$$
Moreover, $\hat{m}_+ = m_+ m_{0+}^{-1}$ and $\hat{m}_- = m_-$ on $\gamma_+$, while $\hat{m}_+ = m_+$ and $\hat{m}_- = m_-m_{0-}^{-1}$ on $\gamma_-$.
It follows that $\hat{m}_+ = \hat{m}_- \hat{v}$ a.e. on $\hat{\Gamma}$.
Hence $\hat{m}$ satisfies the $L^p$-RH problem determined by $(\hat{\Gamma}, \hat{v})$.

Conversely, suppose $\hat{m} \in I + \dot{E}^p(\hat{D})$ satisfies the $L^p$-RH problem determined by $(\hat{\Gamma}, \hat{v})$ and define $m(z)$ for $z \in \C \setminus \hat{\Gamma}$  by (\ref{Amdef}). By Lemma \ref{fglemma}, $m \in I + \dot{E}^p(\hat{D})$.
The nontangential boundary values $m_\pm \in I + \dot{L}^p(\hat{\Gamma})$ satisfy
$$m_\pm = \begin{cases}
\hat{m}_\pm m_{0\pm} & \text{on} \quad \hat{\Gamma} \cap B_+,\\
\hat{m}_\pm, & \text{on} \quad \hat{\Gamma} \cap B_-.
\end{cases}$$
Moreover, $m_+ = \hat{m}_+ m_{0+}$ and $m_- = \hat{m}_-$ on $\gamma_+$, while $m_+ = \hat{m}_+$ and $m_- = \hat{m}_-m_{0-}$ on $\gamma_-$.
It follows that $m_+ = m_- v$ a.e. on $\Gamma$ and that $m_+ = m_-$ a.e. on $\gamma$. 
Using Theorem \ref{jumpth1} and the fact that $m_+ = m_-$ a.e. on $\gamma$, we find
\begin{align}\label{mminusIChat}
m(z) - I = (\mathcal{C}_{\hat{\Gamma}}(m_+ - m_-))(z) = (\mathcal{C}_{\Gamma}(m_+ - m_-))(z), \qquad z \in \hat{D}.
\end{align}
Since the right-hand side belongs to $\dot{E}^p(D)$ by part $(b)$ of Theorem \ref{jumpth1}, defining $m(z)$ for $z \in \gamma$ by $m(z) = m_+(z) = m_-(z)$, we have $m \in I + \dot{E}^p(D)$ and equation (\ref{mminusIChat}) becomes valid for all $z \in D$. It follows that $m$ satisfies the $L^p$-RH problem determined by $(\Gamma, v)$.
\proofend

\section{Conclusions}\nequation
We have developed a theory of $L^p$-Riemann-Hilbert problems for a class of jump contours of low regularity. More precisely, we have considered jump contours $\Gamma$ which are the union of a finite number of simple closed Carleson curves. Several results well-known from the case of smooth contours have been shown to generalize to this more general setting. Our definition of a solution of the $L^p$-RH problem has been novel in that it 
has been given directly in terms of $m(z)$ using appropriate Smirnoff classes (and not in terms of $m_\pm$ as in \cite{DZ2002a, DZ2002b, FIKN2006, Z1989}). Moreover, we have established uniqueness of the $L^p$-RH problem for $n\times n$ matrices for any $1 \leq n \leq p$ (see Theorem \ref{uniqueth}; for $n = p = 2$ this result was proved in \cite{D1999, DZ2003} for sufficiently smooth contours). Overall it has been demonstrated that the theory of $L^p$-RH problems extends virtually unimpeded to the setting of Carleson jump contours.

On the other hand, it is natural to expect the class of Carleson contours to be the largest class of contours for which a clean RH theory exists. Indeed, the Cauchy singular operator $\mathcal{S}_\Gamma$, which is essential in the RH formalism, is known to be bounded on $L^p(\Gamma)$, $1 < p < \infty$, if and only if $\Gamma$ is a Carleson curve \cite{BK1997}.

The presented results can be used to determine rigorously the long-time asymptotics of solutions of integrable evolution equations via the method of nonlinear steepest descent.
We mention in this regard that RH problems with complicated contours that do not fit into the traditional framework arise in the analysis of initial-boundary value problems for integrable PDEs. For example, the analysis of the Degasperis-Procesi equation on the half-line leads to a RH problem with a jump contour involving nontransversal intersections, see Figure \ref{Dns.pdf}.

\appendix
\section{Proof of Proposition \ref{Jprop}} \label{Japp}
\renewcommand{\theequation}{A.\arabic{equation}}\nequation
We first prove a lemma. 

\begin{lemma}\label{applemma1}
Let $\Gamma \subset \C$ be an arc homeomorphic to $I$ where $I$ is either $[0,1]$, $[0,1)$, or $(0,1]$. If $\gamma:I \to \Gamma$ is a homeomorphism, then $\Gamma$ is locally rectifiable if and only if $\gamma((0,1))$ is locally rectifiable. 
\end{lemma}
\proofbegin
We will prove that $\Gamma$ is rectifiable whenever $I = [0,1]$ and $\gamma((0,1))$ is locally rectifiable; the other cases can easily be reduced to this case. Suppose $\gamma((0,1))$ is locally rectifiable. Since $\gamma((0,1))$ is contained in the bounded set $\gamma([0,1])$, $\gamma((0,1))$ is rectifiable. Let $a = t_0 < t_1 < \dots < t_N = b$ be a partition of a closed subinterval $[a,b]$ of $(0,1)$. Since $\gamma((0,1))$ is rectifiable, 
$$\sup_{0 < a < b < 1} \sup_{\substack{\text{partitions} \\ \text{$\{t_i\}$ of $[a,b]$}}} \sum_{j=1}^N |\gamma(t_j) - \gamma(t_{j-1})| < \infty.$$
On the other hand, since $\gamma([0,1])$ is compact, 
$$\sup_{0 < a < b < 1} (|\gamma(a) - \gamma(0)| + |\gamma(1) - \gamma(b)|) < \infty.$$
Thus,
$$\sup_{0 < a < b < 1} \sup_{\substack{\text{partitions} \\ \text{$\{t_i\}$ of $[a,b]$}}} \bigg(|\gamma(a) - \gamma(0)| + |\gamma(1) - \gamma(b)| + \sum_{j=1}^N |\gamma(t_j) - \gamma(t_{j-1})|\bigg) < \infty,$$
showing that $\Gamma$ is rectifiable.
\proofend

We now prove Proposition \ref{Jprop}. Let $\Gamma \subset \hat{\C}$ be a Carleson curve, that is, $\Gamma$ is connected and $\Gamma \cap \C$ Êis a locally rectifiable composed curve satisfying (\ref{carlesondef}). We need to prove that $\psi(\Gamma)  \subset \hat{\C}$ is a Carleson curve. The proof is trivial if $c = 0$. Thus suppose $c \neq 0$.
Since $\psi = \psi_4 \circ \psi_3 \circ \psi_2 \circ \psi_1$ is the composition of the four maps
$$\psi_1(z) = z + \frac{d}{c}, \quad
\psi_2(z) = \frac{1}{z}, \quad
\psi_3(z) = \frac{bc - ad}{c^2} z, \quad
\psi_4(z) = z + \frac{a}{c},$$
and the operations of multiplication and translation by a complex number clearly preserve the family of Carleson curves, we may assume that $\psi(z) = z^{-1}$. Since $\psi(\Gamma) \cap \C$ is Carleson if and only if each of its finite number of arcs is Carleson, we may assume that $\Gamma$ consists of a single (possibly unbounded) arc and that $\Gamma \subset \C$. Furthermore, Lemma \ref{applemma1} shows that we may discard any possible endpoints, and hence assume that $\Gamma$ is homeomorphic to the open interval $(0,1)$. 
Finally, if $0 \in \Gamma$, we may consider each of the two arcs that make up $\Gamma \setminus \{0\}$ separately. 
Thus, without loss of generality, let $\psi(z)  = z^{-1}$ and let $\Gamma \subset \C$ be an arc homeomorphic to $(0,1)$ such that $0 \notin \Gamma$. Then $\psi(\Gamma)\subset \C$ is an arc homeomorphic to $(0,1)$ such that $0 \notin \psi(\Gamma)$. We need to prove that $\psi(\Gamma)$ is locally rectifiable and Carleson. Let $\gamma:(0,1) \to \Gamma$ be a homeomorphism.

\begin{lemma}\label{Gamma0lemma}
If $\Gamma_0$ is a subarc of $\Gamma$ such that there exist constants $m, M \in (0, \infty)$ with the property that $m \leq |z| \leq M$ for all $z \in \Gamma_0$, then both $\Gamma_0$ and $\psi(\Gamma_0)$ are rectifiable and
\begin{align}\label{m2Gamma}
m^2 |\psi(\Gamma_0)| \leq |\Gamma_0| \leq M^2 |\psi(\Gamma_0)|.
\end{align}
\end{lemma}
\proofbegin
Since $\Gamma$ is locally rectifiable and $\Gamma_0$ is a bounded subarc, $\Gamma_0$ is rectifiable. 
Let $I \subset (0,1)$ be the subinterval of $(0,1)$ for which $\Gamma_0 = \gamma(I)$. Then $I$ equals $[a,b]$, $[a,b)$, $(a,b]$, or $(a,b)$ for some $0< a < b < 1$.
If $c \leq t_0 < t_1 < \cdots < t_{N} \leq d$ is a partition of a closed subinterval $[c,d] \subset I$, then
$$\sum_{j=1}^N |\gamma(t_j) - \gamma(t_{j-1})|
= \sum_{j=1}^N |\psi(\gamma(t_j)) - \psi(\gamma(t_{j-1}))| |\gamma(t_j)| |\gamma(t_{j-1})|$$
and so
$$m^2 \sum_{j=1}^N |\psi(\gamma(t_j)) - \psi(\gamma(t_{j-1}))| 
\leq \sum_{j=1}^N |\gamma(t_j) - \gamma(t_{j-1})|
\leq M^2 \sum_{j=1}^N |\psi(\gamma(t_j)) - \psi(\gamma(t_{j-1}))|.$$
Taking the supremum over all partitions and all closed subintervals $[c,d] \subset I$, we find (\ref{m2Gamma}). 
\proofend

We next prove that $\psi(\Gamma)$ is locally rectifiable. Let $[a,b]$ be a closed subinterval of $(0,1)$. Let $\Gamma_c = \gamma([a,b])$. Since $\Gamma_c$ is compact and $0 \notin \Gamma_c$, $\Gamma_c$ is bounded and bounded away from $0$; hence $\psi(\Gamma_c)$ is also bounded and bounded and bounded away from $0$. It follows that $\psi(\Gamma_c)$ is rectifiable and Lemma \ref{Gamma0lemma} implies
\begin{align} \nonumber
|\psi(\Gamma_c) \cap D(0, r)| & = \sum_{n=1}^\infty \big|\psi(\Gamma_c) \cap \{2^{-n}r \leq |w| < 2^{1-n}r\}\big|
	\\ \label{psiCestimate}
& \leq \sum_{n=1}^\infty r^2 2^{2 - 2n} \big|\Gamma_c \cap \{z|2^{n-1} r^{-1} < |z| \leq 2^{n} r^{-1}\} \big|.
\end{align}
By the Carleson property (\ref{carlesondef2}) of $\Gamma$ applied to $D(0, 2^{n} r^{-1})$, the right-hand side of (\ref{psiCestimate}) is bounded above by
$$C_\Gamma \sum_{n=1}^\infty r^2 2^{2 - 2n} 2^n r^{-1} = 4 C_\Gamma r,$$
where $C_\Gamma > 0$ is a constant. 
Since the closed interval $[a,b] \subset (0,1)$ was arbitrary, it follows that 
$$\sup_{0< a < b < 1} |\psi(\gamma([a,b])) \cap D(0, r)| \leq 4 C_\Gamma r < \infty$$
for each $r>0$. This shows that $\psi(\Gamma)$ is locally rectifiable.

It remains to prove that $\psi(\Gamma)$ is Carleson. Let $w_0 \in \psi(\Gamma)$ and $r > 0$. Let $R = |w_0| + r$. Then
\begin{align*}
|\psi(\Gamma) \cap D(w_0, r)|
& = |\psi(\Gamma) \cap D(w_0, r) \cap D(0, R) |
	\\
& = \sum_{n=1}^\infty \big|\psi(\Gamma) \cap D(w_0, r) \cap \{w| 2^{-n}R \leq |w| < 2^{1-n} R\} \big|.
\end{align*}
In view of (\ref{m2Gamma}), this yields
$$|\psi(\Gamma) \cap D(w_0, r)|
\leq \sum_{n=1}^\infty 2^{2-2n}R^2 \big|\Gamma \cap \{z||z^{-1} - z_0^{-1}| < r\} \cap \{z | 2^{n-1}R^{-1} < |z| \leq 2^{n} R^{-1}\} \big|,$$
where $z_0 = w_0^{-1}$.
The set $\Gamma \cap \{z | |z^{-1} - z_0^{-1}| < r\} \cap \{z | 2^{n-1}R^{-1} < |z| \leq 2^{n} R^{-1}\}$ is contained in the intersection of the two open disks $D(z_0, r |z_0| 2^{n} R^{-1})$ and $D(0, 2^{n+1} R^{-1})$.
Hence we may use the Carleson property of $\Gamma$ on these disks to find
\begin{align*}
|\psi(\Gamma) \cap D(w_0, r)|
& \leq C_\Gamma \sum_{n=1}^\infty 2^{2-2n}R^2 \min\big\{ r |z_0| 2^{n} R^{-1}, 2^{n+1} R^{-1}\big\}
	\\
& \leq 4 C_\Gamma R \min\{r |z_0|, 2\}
= 4 C_\Gamma \min\{r(1 + r |w_0|^{-1}), 2(|w_0| + r)\}.
\end{align*}
where $C_\Gamma > 0$ is a constant. 
If $|w_0| \geq r$, then $4 C_\Gamma r(1 + r |w_0|^{-1}) \leq 8 C_\Gamma r$.
If $|w_0| < r$, then $8 C_\Gamma (|w_0| + r) < 16 C_\Gamma r$. Hence
$$|\psi(\Gamma) \cap D(w_0, r)| < 16 C_\Gamma r,$$
for all $w_0 \in \psi(\Gamma)$ and all $r > 0$. This proves that $\psi(\Gamma)$ is Carleson and completes the proof of Proposition \ref{Jprop}.
\proofend

\bigskip
\noindent
{\bf Acknowledgement} {\it The author is grateful to Prof. I. M. Spitkovsky for helpful remarks on a first version of the manuscript and acknowledges support from the EPSRC, UK.}

\bibliographystyle{plain}
\bibliography{is}

\begin{thebibliography}{99}
\small

\bibitem{AF2003}
M. J. Ablowitz and A. S. Fokas, {\it Complex variables: introduction and applications.} 2nd edition. Cambridge Texts in Applied Mathematics. Cambridge University Press, Cambridge, 2003.

\bibitem{BC1984}
R. Beals and R. R. Coifman, Scattering and inverse scattering for first order systems, {\it Comm. Pure Appl. Math.} {\bf 37} (1984), 39--90.

\bibitem{BK1997}
A. B\"ottcher and Y. I. Karlovich, {\it Carleson curves, Muckenhoupt weights, and Toeplitz operators}, Progress in Mathematics, 154. Birkh\"auser Verlag, Basel, 1997.

\bibitem{CG1981}
K. F. Clancey and I. Gohberg, {\it Factorization of matrix functions and singular integral operators.} Operator Theory: Advances and Applications, 3. BirkhŠuser Verlag, Basel-Boston, Mass., 1981.

\bibitem{D1984}
G. David, Op\'erateurs int\'egraux singuliers sur certaines courbes du plan complexe, {\it Ann. Sci. \'Ecole Norm. Sup. (4)} {\bf 17} (1984), 157--189.

\bibitem{D1999}
P. Deift, {\it Orthogonal polynomials and random matrices: a Riemann-Hilbert approach}, Courant Lecture Notes in Mathematics, 3, New York University, Courant Institute of Mathematical Sciences, New York; American Mathematical Society, Providence, RI, 1999.
 
\bibitem{DVZ1997}
P. Deift, S. Venakides, and X. Zhou, New results in small dispersion KdV by an extension of the steepest descent method for Riemann-Hilbert problems, {\it Internat. Math. Res. Notices} {\bf 1997}, 286--299. 

\bibitem{DZ1993}
P. Deift and X. Zhou, A steepest descent method for oscillatory Riemann-
Hilbert problems. Asymptotics for the MKdV equation, 
{\it Ann. of Math.} {\bf 137} (1993), 295--368.

\bibitem{DZ2002a}
P. Deift and X. Zhou, Perturbation theory for infinite-dimensional integrable systems on the line. A case study, {\it Acta Math.} {\bf 188} (2002), 163--262.

\bibitem{DZ2002b}
P. Deift and X. Zhou, A priori $L^p$-estimates for solutions of Riemann-Hilbert problems, {\it Int. Math. Res. Not.} {\bf 2002}, 2121--2154.

\bibitem{DZ2003}
P. Deift and X. Zhou, Long-time asymptotics for solutions of the NLS equation with initial data in a weighted Sobolev space, {\it Comm. Pure Appl. Math.} {\bf 2003}, 1029--1077.

\bibitem{D1970}
P. L. Duren, {\it Theory of $H^p$ spaces}, Pure and Applied Mathematics, Vol. 38 Academic Press, 1970.

\bibitem{FI1996}
A. S. Fokas and A. R. Its, The linearization of the initial-boundary value problem of the nonlinear Schr\"odinger equation,
{\it SIAM J. Math. Anal.} {\bf 27} (1996), 738--764. 

\bibitem{FIKN2006}
A. S. Fokas, A. R. Its, A. A. Kapaev, V. Y. Novokshenov, {\it Painlev\'e transcendents. The Riemann-Hilbert approach.} Mathematical Surveys and Monographs, 128. American Mathematical Society, Providence, RI, 2006.

\bibitem{G1987}
D. Gaier, {\it Lectures on complex approximation.} Birkh\"auser Boston, Boston, 1987.

\bibitem{GKS2003}
I. Gohberg, M. A. Kaashoek, and I. M. Spitkovsky, An overview of matrix factorization theory and operator applications, Factorization and integrable systems (Faro, 2000), 1--102, 
Oper. Theory Adv. Appl., 141, Birkh\"auser, Basel, 2003.

\bibitem{I1981}
A. R. Its, Asymptotic behavior of the solutions to the nonlinear Schr\"odinger equation, and isomonodromic deformations of systems of linear differential equations, {\it Dokl. Akad. Nauk SSSR} {\bf 261} (1981), 14--18 (in Russian); {\it Soviet Math. Dokl.} {\bf 24} (1982), 452--456 (in English).

\bibitem{K2008}
S. Kamvissis, From stationary phase to steepest descent. Integrable systems and random matrices, 145--162, 
Contemp. Math., 458, Amer. Math. Soc., Providence, RI, 2008. 

\bibitem{KT2012}
S. Kamvissis and G. Teschl, Long-time asymptotics of the periodic Toda lattice under short-range perturbations, {\it J. Math. Phys.} {\bf 53} (2012), 073706, 35 pp.

\bibitem{KKP1998} 
G. Khuskivadze, V. Kokilashvili, and V. Paatashvili, Boundary value problems for analytic and harmonic functions in domains with nonsmooth boundaries. Applications to conformal mappings. {\it Mem. Differential Equations Math. Phys.} {\bf 14} (1998), 195 pp.
 
\bibitem{L2013}
J. Lenells, The Degasperis-Procesi equation on the half-line, {\it Nonlinear Anal.} {\bf 76} (2013), 122--139.
 
\bibitem{LS1987}
G. S. Litvinchuk and I. M. Spitkovskii, {\it Factorization of measurable matrix functions}, Operator Theory: Advances and Applications 25, Birkh\"auser Verlag, Basel, 1987.

\bibitem{M1974}
S. V. Manakov, Nonlinear Fraunhofer diffraction, {\it Zh. Eksp. Teor. Fiz.} {\bf 65} (1973), 1392--1398 (in Russian); {\it Sov. Phys. JETP} {\bf 38} (1974), 693--696 (in English).

\bibitem{M1992}
N. I. Muskhelishvili, {\it Singular integral equations. Boundary problems of function theory and their application to mathematical physics.} Dover Publications, New York, 1992.
 
\bibitem{P1956}
I. I. Priwalow, {\it Randeigenschaften analytischer Funktionen}, Hochschulb\"ucher f\"ur Mathematik, Bd. 25. VEB Deutscher Verlag der Wissenschaften, Berlin, 1956.

\bibitem{R1988}
Y. L. Rodin, {\it The Riemann boundary problem on Riemann surfaces.} Mathematics and its Applications (Soviet Series), 16. D. Reidel Publishing Co., Dordrecht, 1988.

\bibitem{Z1989}
X. Zhou, The Riemann-Hilbert problem and inverse scattering, {\it SIAM J. Math. Anal.} {\bf 20} (1989), 966--986. 

\bibitem{Z1971}
\`E. I. Zverovi\v{c}, Boundary value problems in the theory of analytic functions in H\"older classes on Riemann surfaces, Russian Math. Surveys \bf{26} (1971), 117--192.

\end{thebibliography}

\end{document}